\pdfoutput=1
\RequirePackage{ifpdf}
\ifpdf 
\documentclass[pdftex]{sigma}
\else
\documentclass{sigma}
\fi

\numberwithin{equation}{section}

\newtheorem{Theorem}{Theorem}[section]

\newtheorem{Lemma}[Theorem]{Lemma}
\newtheorem{Proposition}[Theorem]{Proposition}
 { \theoremstyle{definition}

\newtheorem{Remark}[Theorem]{Remark} }

\usepackage[labelsep=period,labelfont=bf]{caption}
\usepackage{subcaption}
\usepackage{bbm}
\usepackage{bm}
\usepackage{mathtools}
\usepackage{enumitem}
\usepackage{tikz-cd}

\newcommand{\RN}[1]{%
	\textup{\uppercase\expandafter{\romannumeral#1}}%
}

\def\bfs{\boldsymbol}

\def\pa{\partial}

\def\wh{\widehat}
\def\wt{\widetilde}

\def\C{\mathbb{C}}

\def\P{\mathbf{P}}
\def\G{\mathbf{G}}
\def\R{\mathbb{R}}

\newcommand{\Pf}{{\textup{Pf}}}
\newcommand{\erfc}{\operatorname{erfc}}
\newcommand{\erf}{\operatorname{erf}}
\newcommand{\bfR}{\mathbf{R}}

\newcommand{\bfkappa}{{\bm \varkappa}}

\newcommand{\Prob}{{\mathbb{P}}}

\begin{document}
\allowdisplaybreaks

\renewcommand{\thefootnote}{}

\newcommand{\arXivNumber}{2209.01934}

\renewcommand{\PaperNumber}{033}

\FirstPageHeading

\ShortArticleName{Spherical Induced Ensembles with Symplectic Symmetry}

\ArticleName{Spherical Induced Ensembles \\ with Symplectic Symmetry\footnote{This paper is a~contribution to the Special Issue on Evolution Equations, Exactly Solvable Models and Random Matrices in honor of Alexander Its' 70th birthday. The~full collection is available at \href{https://www.emis.de/journals/SIGMA/Its.html}{https://www.emis.de/journals/SIGMA/Its.html}}}

\Author{Sung-Soo BYUN~$^{\rm a}$ and Peter J. FORRESTER~$^{\rm b}$}

\AuthorNameForHeading{S.-S.~Byun and P.J.~Forrester}

\Address{$^{\rm a)}$~Center for Mathematical Challenges, Korea Institute for Advanced Study,\\
\hphantom{$^{\rm a)}$}~Seoul 02455, Republic of Korea}
\EmailD{\href{mailto:email@address}{sungsoobyun@kias.re.kr}}

\Address{$^{\rm b)}$~School of Mathematical and Statistics, The University of Melbourne, Victoria 3010, Australia}
\EmailD{\href{mailto:email@address}{pjforr@unimelb.edu.au}}

\ArticleDates{Received September 22, 2022, in final form May 16, 2023; Published online May 30, 2023}

\Abstract{We consider the complex eigenvalues of the induced spherical Ginibre ensemble with symplectic symmetry and establish the local universality of these point processes along the real axis. We derive scaling limits of all correlation functions at regular points both in the strong and weak non-unitary regimes as well as at the origin having spectral singularity. A~key ingredient of our proof is a derivation of a differential equation satisfied by the correlation kernels of the associated Pfaffian point processes, thereby allowing us to perform asymptotic analysis.}


\Classification{60B20; 33C45; 33E12}

 \begin{flushright}
 \begin{minipage}{60mm}
 \textit{Dedicated to Alexander Its\\ on the occasion of his 70th birthday}
 \end{minipage}
 \end{flushright}

\renewcommand{\thefootnote}{\arabic{footnote}}
\setcounter{footnote}{0}

\section{Introduction and results}

The footprints of the universality in non-Hermitian random matrix theory began with the work~\cite{ginibre1965statistical} of Ginibre.
There, three classes of Gaussian random matrices with complex, real, and quaternion elements were introduced, and they are now called the Ginibre ensembles.
(We refer to~\cite{byun2022progress,byun2023progress} for recent reviews on these topics.)
Although the eigenvalues of the matrices in each symmetry class all follow the universal circular law at the macroscopic level, their statistical properties are quite different from many perspectives.
For instance, in the complex symmetry class, the real axis is not special due to the rotational invariance.
On the other hand, in the real and quaternion cases, there exist microscopic attraction and repulsion respectively along the real axis.\looseness=-1

The difference among these three symmetry classes can also be found in their integrable structures.
More precisely, the eigenvalues of the complex matrices form determinantal point processes, whereas those of the real and quaternion matrices form Pfaffian point processes.
Furthermore, while the correlation kernels of the complex matrices can be written in terms of the planar orthogonal polynomials, their counterparts for the real and quaternion matrices are described in terms of the (less understood) planar skew-orthogonal polynomials.

Due to the more complicated integrable structures of Pfaffian point processes, it is not surprising that the local universality classes (i.e., scaling limits of all eigenvalue correlation functions) were first investigated in the complex symmetry class.
Indeed, the bulk scaling limit of the complex Ginibre ensemble was already introduced in the work \cite{ginibre1965statistical} of Ginibre.
On the other hand, the edge scaling limit of the complex Ginibre ensemble was discovered in \cite{MR1690355}.
For the real symmetry class, the bulk and edge scaling limits of the Ginibre ensemble were investigated in \cite{MR2530159,MR2430570, MR2371225}.
Finally, for the quaternion case, the bulk scaling limit was first introduced in the second edition of Mehta's book~\cite{Mehta} and later rediscovered by Kanzieper~\cite{MR1928853}.
In contrast, the edge scaling limit in this symmetry class was discovered only recently in~\cite{akemann2021scaling}.
(See also~\cite{Lysychkin} for an alternative derivation for the $1$-point function.)

From the above philosophy, it is not surprising again that the universality principle was first established in the complex symmetry class.
Among plenty of works in this direction, the bulk universality of random normal matrix ensembles was obtained in \cite{MR2817648}.
More recently, the edge universality of these models was obtained in \cite{hedenmalm2017planar}, where the authors developed a~general asymptotic theory of the planar orthogonal polynomials.
However, the literature on the universality in the other symmetry classes are more limited.

Nevertheless, there have been several recent works on the scaling limits of \emph{planar symplectic ensembles}, which are contained in the symmetry class of the quaternion Ginibre ensemble. (By definition, these are point processes which follow the joint probability distribution~\eqref{Gibbs}.)
For instance, the universal scaling limits of the symplectic elliptic Ginibre ensemble at the origin were obtained in~\cite{akemann2021skew} and were extended in~\cite{byun2021universal} along the whole real axis.
Furthermore, non-standard universality classes under the presence of certain singularities have been discovered as well.
To name a few, the scaling limits at the \emph{singular origin} were studied in~\cite{akemann2021scaling} for the Mittag-Leffler ensemble (a generalisation of the symplectic induced Ginibre ensemble), in~\cite{MR3066113} for the product ensembles and in~\cite{MR2180006,akemann2021skew} for the Laguerre ensembles.
The boundary scaling limits under the \emph{hard edge} type conditions were investigated in \cite{khoruzhenko2021truncations} for the truncated symplectic unitary ensembles and in \cite{byun2021wronskian} for the Ginbire ensemble with boundary confinements.
Beyond the above-mentioned cases, the scaling limits of the models interpolating one- and two-dimensional ensembles have also been studied.
In this direction, the scaling limits of the symplectic elliptic Ginibre ensemble in the \emph{almost-Hermitian} (or \emph{weakly non-Hermitian}) regime were derived in \cite{MR3192169,byun2021wronskian, MR1928853}.
Very recently, the scaling limits of the symplectic induced Ginibre ensemble in the \emph{almost-circular} (or \emph{weakly non-unitary}) regime were obtained in \cite{byun2022almost}.
While the almost-Hermitian \cite{akemann2016universality,fyodorov1997almost, MR1634312} and almost-circular \cite{AB21,byun2021random} ensembles have the same bulk scaling limits in the complex symmetry class, those are different in the symplectic symmetry class in the vicinity of the real line due to the lack of the translation invariance; see \cite{byun2022almost} for further details.

\begin{figure}[h!]
		\begin{subfigure}{0.32\textwidth}
		\begin{center}	
			\includegraphics[width=\textwidth]{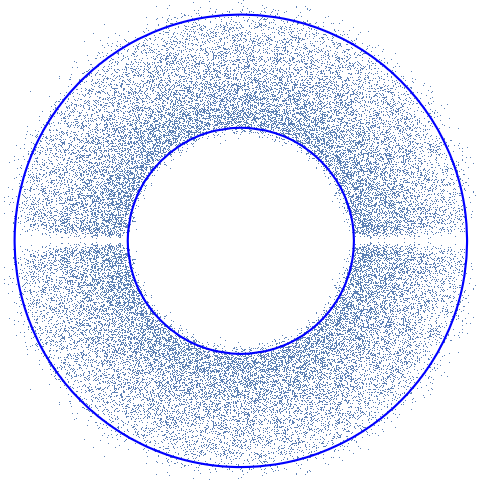}
		\end{center}
		\subcaption{$L=N$, $n=2N$
		}
	\end{subfigure}	
		\begin{subfigure}{0.33\textwidth}
		\begin{center}	
			\includegraphics[width=\textwidth]{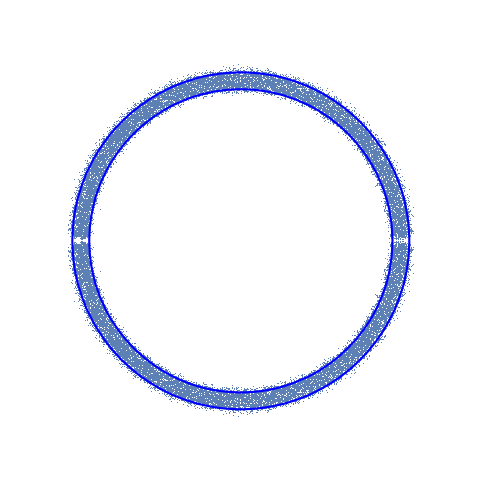}
		\end{center}
		\subcaption{$L=\frac{N^2}{\rho^2}-N$, $n=\frac{N^2}{\rho^2}$, $\rho=\sqrt{10}$}
	\end{subfigure}	
		\begin{subfigure}{0.32\textwidth}
		\begin{center}	
			\includegraphics[width=\textwidth]{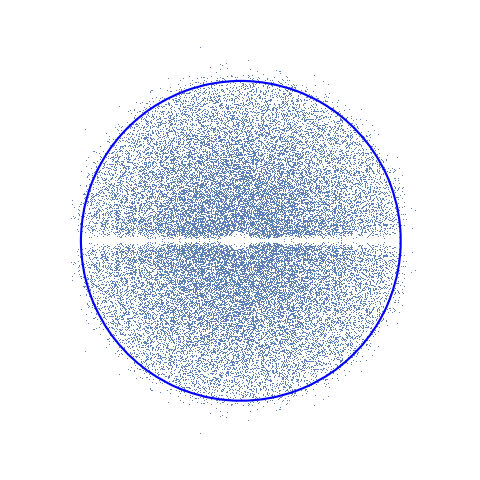}
		\end{center}
		\subcaption{$L=1$, $n=2N$}
	\end{subfigure}	

	\begin{subfigure}{0.32\textwidth}
		\begin{center}	
			\includegraphics[height=\textwidth]{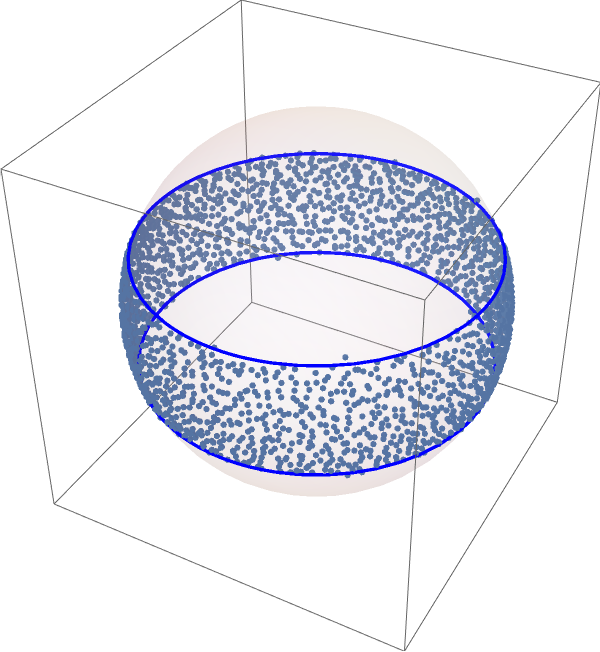}
		\end{center}
		\subcaption{$L=N$, $n=2N$}
	\end{subfigure}	
		\begin{subfigure}{0.33\textwidth}
		\begin{center}	
			\includegraphics[height=\textwidth]{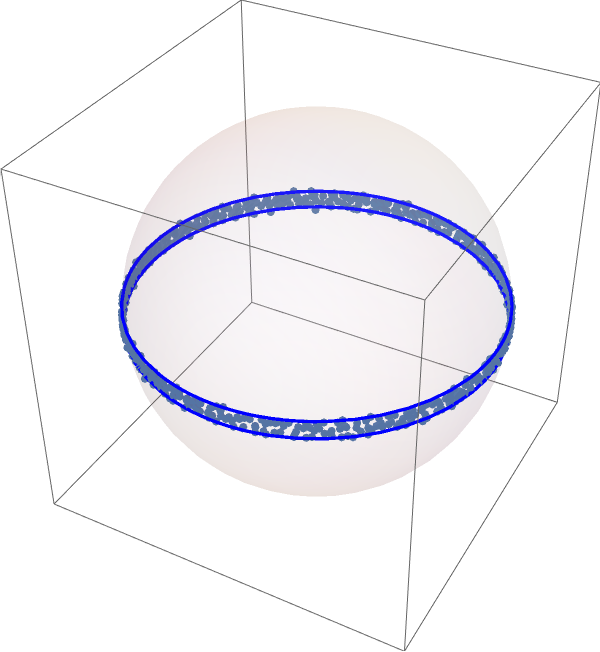}
		\end{center}
		\subcaption{$L=\frac{N^2}{\rho^2}-N$, $n=\frac{N^2}{\rho^2}$, $\rho=\sqrt{10}$}
	\end{subfigure}	
		\begin{subfigure}{0.32\textwidth}
		\begin{center}	
			\includegraphics[height=\textwidth]{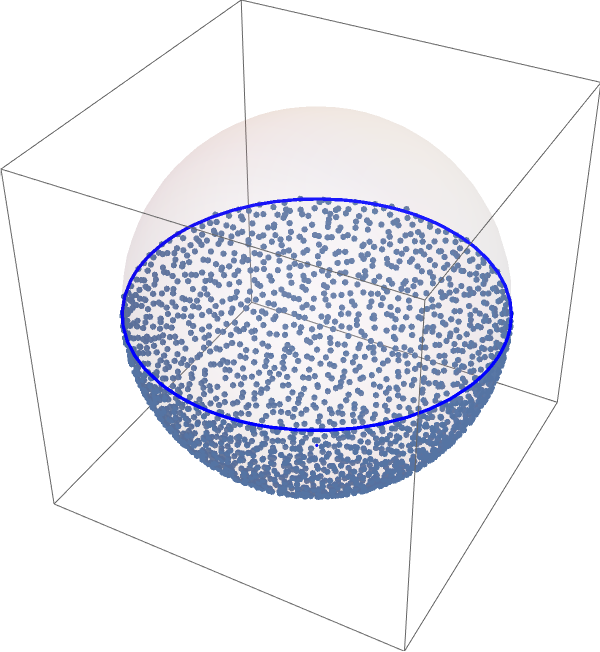}
		\end{center}
		\subcaption{$L=1$, $n=2N$}
	\end{subfigure}	
	\caption{The plots (a)--(c) display the eigenvalues of $\textbf{G}$ for $N=100$ and $200$ realisations with different values of $n$ and $L$. In all figures (a)--(c), the local repulsion along the real axis is visible.	The plots (d)--(f) show a sample of the eigenvalues of $\textbf{G}$ for $N=1000$ projected onto the unit sphere.} \label{Fig_eigenvalues}
\end{figure}

In this work, we study the \textit{symplectic induced spherical ensembles} with the goal to derive their scaling limits in various regimes and to establish the universality of these point processes.
The symplectic induced spherical ensemble $\G$ is an $N \times N$ quaternion matrix, which is defined by the matrix probability distribution function proportional to
\begin{equation} \label{matrix model}
 \frac{ \det \big(\G \G^\dagger \big)^{ 2L } }{ \det \big( \mathbf{1}_N+\G \G^\dagger \big)^{ 2(n+N+L) } }.
\end{equation}
 Here $n$ and $L$ are parameters, with $n \ge N$ and $L \ge 0 $ also possibly dependent on $N$.
In particular, if $n=N$, $L=0$, the model \eqref{matrix model} is known as the spherical ensemble with symplectic symmetry.
The name ``spherical'' originates from the fact that their eigenvalues tend to be uniformly distributed on the unit sphere under the (inverse) stereographic projection; see, e.g., \cite{forrester2012pfaffian,Krishnapur09}. And as discussed in the ensuing text, the term symplectic symmetry relates to an invariance of the underlying Gaussian matrices.
To realise the matrix probability distribution~\eqref{matrix model}, following~\cite{MR2881072} and~\cite[Appendix~B]{MR3612266} first introduce a particular $(N + L) \times N$ random matrix~$\mathbf Y$ with each entry itself a $2 \times 2$ matrix representation of a quaternion.
The matrix is said to have quaternion entries for short; see, e.g.,~\cite[Section 1.3.2]{forrester2010log}.
The specification of $\mathbf Y$ is that $\mathbf Y = \mathbf X \mathbf A^{-1/2}$, where $\mathbf X$ is an $(N+L) \times N$ standard
Gaussian matrix with quaternion entries (also referred to as a rectangular quaternion Ginibre matrix), while $\mathbf A$ is an $N \times N$ Wishart matrix with quaternion entries. More explicitly, $\mathbf A = \mathbf Q^\dagger \mathbf Q$ where~$\mathbf Q$ is an $n \times N$ rectangular quaternion Ginibre matrix; see, e.g.,~\cite[Definition~3.1.2]{forrester2010log}.
In terms of such $\mathbf Y$, and a Haar distributed unitary random matrix with quaternion entries (i.e., a symplectic unitary matrix $\mathbf U$; see, e.g.,~\cite{DF17}), define $\mathbf G = \mathbf U (\mathbf Y^\dagger \mathbf Y)^{1/2}$.
It is the random matrix $\mathbf G$ which has matrix distribution~\eqref{matrix model}; the corresponding eigenvalues, which must come in complex conjugate pairs, are used in producing the plots of Figure~\ref{Fig_eigenvalues}.

Several fundamental properties of the symplectic induced spherical ensemble were discovered by Mays and Ponsaing \cite{MR3612266}.
(We also refer to an earlier work \cite{mays2013real} on the induced spherical ensemble with orthogonal symmetry.)
In particular, it was shown in \cite[Section 3]{MR3612266} that the joint probability distribution $\P_N$ of its independent eigenvalues $\bfs{\zeta}= \{ \zeta_j \}_{j=1}^N$ is given by
\begin{gather}\label{Gibbs}
{\rm d}\P_N(\boldsymbol{\zeta}) = \frac{1}{N!Z_N} \prod_{1 \leq j<k \leq N} |\zeta_j-\zeta_k|^2 \big|\zeta_j-\overline{\zeta}_k\big|^2 \prod_{j=1}^{N} \big|\zeta_j-\overline{\zeta}_j\big|^2 {\rm e}^{ -2N Q(\zeta_j) } \, {\rm d}A(\zeta_j),
\end{gather}
where ${\rm d}A(\zeta):={\rm d}^2\zeta/\pi$, and $Z_N$ is the normalisation constant.
While the independent eigenvalues should have each $\zeta_j$ in the upper half complex plane, relaxing this condition leaves~\eqref{Gibbs} unaltered and simplifies the presentation.
Here the potential $Q$ is given by
\begin{equation} \label{Q potential}
Q(\zeta):= \frac{n+L+1}{N}\log\big(1+|\zeta|^{2}\big)-\frac{2L}{N}\log |\zeta|.
\end{equation}
We remark that the distribution \eqref{Gibbs} can be interpreted as a two-dimensional Coulomb gas ensemble \cite{forrester2016analogies,kiessling1999note} with additional complex conjugation symmetry; see also Appendix~\ref{A1}.

We first briefly recall the macroscopic property of the ensemble \eqref{Gibbs}.
Combining the convergence of the empirical measure \cite{MR2934715} and the basic facts from the logarithmic potential theory \cite[Section~IV.6]{ST97}, one can see that as $N \to \infty$, the eigenvalues $\bfs{\zeta}$ tend to be distributed on the droplet
\begin{equation} \label{droplet}
S=\{\zeta \in \C\colon r_1 \le |\zeta| \le r_2\},\qquad r_1= \sqrt{ \frac{L}{n} }, \qquad r_2= \sqrt{ \frac{N+L}{n-N} }
\end{equation}
with the density
\begin{equation} \label{density}
\frac{n+L}{N} \frac{1}{\big(1+|\zeta|^2\big)^2}.
\end{equation}
This property was also shown in \cite[Section~6]{MR3612266} using a different method.
We also refer to~\cite{criado2022vector} and references therein for recent works on the equilibrium measure problems on the sphere under the insertion of point charges.

For detailed statistical information about the ensemble~\eqref{Gibbs}, we study its $k$-point correlation function
\begin{equation}\label{bfRNk def}
	\bfR_{N,k}(\zeta_1,\dots, \zeta_k) := \frac{N!}{(N-k)!} \int_{\C^{N-k}} \P _N(\bfs{\zeta}) \prod_{j=k+1}^N {\rm d}A(\zeta_j).
\end{equation}
The following proposition provides useful formulas to analyse the correlation functions \eqref{bfRNk def}.

\begin{Proposition}[analysis at finite-$N$] \label{Prop_finite N}
For any $N, L \ge 0$, $n \ge N$, and $k \in \mathbb{N}$, the following holds.
\begin{itemize}
 \item[$(a)$] Eigenvalue correlation functions at finite-$N$. We have
\begin{equation}
\label{bfRNk bfkappa tilde}
\bfR_{N,k}(\zeta_1,\dots, \zeta_k) = \Pf \Bigg[
\omega(\zeta_j) \omega(\zeta_l)
	\begin{pmatrix}
		\wt{\bfkappa}_N(\zeta_j,\zeta_l) &
 \wt{\bfkappa}_N\big(\zeta_j,\bar{\zeta}_l\big)
		\\
		\wt{\bfkappa}_N(\bar{\zeta}_j,\zeta_l) & \wt{\bfkappa}_N\big(\bar{\zeta}_j,\bar{\zeta}_l\big)
	\end{pmatrix} \Bigg]_{ j,l=1 }^k \prod_{j=1}^{k} \big(\bar{\zeta}_j-\zeta_j\big),
\end{equation}
where
\begin{equation}
\label{omega}
\omega(\zeta)= \frac{ \big|1+\zeta^2\big|^{n+L-\frac12} }{ 	\big(1+|\zeta|^2\big)^{n+L+1} }.
\end{equation}
Here, the skew-kernel $\wt{\bfkappa}_N(\zeta,\eta)$ is given by
\begin{align} \label{bfkappa tilda G wh}
\wt{\bfkappa}_N(\zeta,\eta) & = \frac{1}{\big(\big(1+\zeta^2\big)\big(1+\eta^2\big)\big)^{n+L-\frac12} } \big( \wh{\boldsymbol{G}}_N(\zeta,\eta)-\wh{\boldsymbol{G}}_N(\eta,\zeta) \big),
\end{align}
where
\begin{equation} \label{GN hat}
\begin{split}
	\wh{\boldsymbol{G}}_N(\zeta,\eta)& := 	\pi \frac{ \Gamma(2n+2L+2) }{ 2^{2L+2n+1} }
	\\
	&\quad \times \sum_{k=0}^{N-1} \sum_{l=0}^{k} \frac{ \zeta^{2k+2L+1} \eta^{2l+2L} }{ \Gamma\big(k+L+\frac32\big) \Gamma(n-k) \Gamma\big(n-l+\frac12\big) \Gamma(l+L+1) } .
\end{split}
\end{equation}
\item[$(b)$] Differential equation for the pre-kernel. The skew-kernel $\wt{\bfkappa}_N(\zeta,\eta)$ satisfies
\begin{equation} \label{CDI}
 \pa_\zeta \wt{\bfkappa}_N(\zeta,\eta) = \frac{1}{\big(1+\zeta^2\big)^{n+L+\frac12}}\big( \RN{1}_N(\zeta,\eta)-\RN{2}_N(\zeta,\eta)-\RN{3}_N(\zeta,\eta) \big),
\end{equation}
where
\begin{align*}
\RN{1}_N(\zeta,\eta) & := \frac{(1+\zeta \eta)^{2n+2L-1} }{\big(1+\eta^2\big)^{n+L-\frac12} } (2n+2L+1)(n+L)
\\
&\quad \times \sum_{k=0}^{2N-1} \binom{2n+2L-1}{k+2L} \mathfrak{p}^{k+2L} (1-\mathfrak{p})^{2n-k-1},
\\
\RN{2}_N(\zeta,\eta) &:= \frac{\zeta^{2N+2L}}{2^{2L+2n}} \frac{ 	 \pi \, \Gamma(2n+2L+2) }{ \Gamma\big(N+L+\frac12\big) \Gamma(n-N) \Gamma\big(n+L+\frac12\big) }
\\
&\quad \times \sum_{k=0}^{N-1} \binom{n+L-\frac12}{k+L} \mathfrak{q}^{k+L} (1-\mathfrak{q})^{n-k-\frac12},
\\
\RN{3}_N(\zeta,\eta) & := \frac{\zeta^{2L-1}}{2^{2L+2n}} \frac{ 	 \pi \, \Gamma(2n+2L+2) }{ \Gamma\big(n+\frac12\big) \Gamma(L) \Gamma\big(n+L+\frac12\big) }
\\
&\quad \times \sum_{k=0}^{N-1} \binom{n+L-\frac12}{k+L+\frac12} \mathfrak{q}^{k+L+\frac12} (1-\mathfrak{q})^{ n-k-1 }.
\end{align*}
Here
\begin{equation} \label{p q mathfrak}
\mathfrak{p}:=\frac{\zeta \eta}{1+\zeta \eta}, \qquad \mathfrak{q}:=\frac{\eta^2}{1+\eta^2}.
\end{equation}
\end{itemize}
\end{Proposition}

\begin{Remark}We stress that Proposition~\ref{Prop_finite N}\,(a) is a direct consequence of the general theory of planar symplectic ensembles~\cite{MR1928853} and the explicit formula of skew-orthogonal polynomials associated with the potential \eqref{Q potential} that can be found in \cite[Proposition~4]{MR3066630}.
(Cf.\ see~\cite[p.~7]{MR3066113} and \cite[Corollary~3.2]{akemann2021skew} for a construction of skew-orthogonal polynomials associated with general radially symmetric potentials.)

Nevertheless, the crux of Proposition~\ref{Prop_finite N} is the transforms \eqref{omega} and \eqref{bfkappa tilda G wh} in the expression~\eqref{bfRNk bfkappa tilde}, which lead to a simple differential equation \eqref{CDI} stated in Proposition~\ref{Prop_finite N}\,(b).
To be more concrete, let us mention that in general, one strategy for performing an asymptotic analysis on a double summation appearing in the skew-orthogonal polynomial kernel is to derive a ``proper'' differential equation satisfied by the kernel; see \cite{MR2180006, akemann2021scaling,byun2022almost, byun2021universal, MR1928853}.
(Such a differential equation for the two-dimensional ensemble is broadly called the generalised Christoffel--Darboux formula \cite{akemann2021scaling,byun2021universal,MR3450566}).
However, if we do not take well-chosen transforms, the resulting differential equation may be difficult to analyse, cf.\ see \cite[Section 6.2]{MR3612266} for a similar discussion on the spherical induced ensemble.

We also mention that the inhomogeneous term $\RN{1}_N(\zeta,\eta)$ in \eqref{CDI} corresponds to the kernel of the complex counterpart \cite{FF11}.
Such a relation has been observed not only for the two-dimensional ensembles \cite{MR2180006,akemann2021scaling,byun2022almost,byun2021universal} but also for their one-dimensional counterparts \cite{MR1762659,MR1675356}.
For a comprehensive summary of this relation for planar ensembles, we refer the reader to \cite{byun2023progress}.
\end{Remark}

\begin{Remark}The terms on the right-hand side of \eqref{CDI} are indeed expressed in a way that one can easily derive their asymptotic behaviours.
	More precisely, the summations in these terms can be written in terms of the incomplete beta functions (see \eqref{RN1 wt 2 beta}, \eqref{RN2 wt 2 beta} and \eqref{RN3 wt 2 beta}) whose asymptotic behaviours are well understood.
	This fact will play an important role in the proof of Theorem~\ref{Thm_main} below.
\end{Remark}

Let us now introduce our main results on various scaling limits of the induced spherical ensembles.
From the microscopic point of view, we first mention that the origin is special since there exists an insertion of a point charge (i.e., $\frac{4L}{N}\log|\zeta|$ term in \eqref{Q potential}, or equivalently the charge $Nq_1$ at the north pole in the sphere picture as given by \eqref{A.3a}), which is also known as the \emph{spectral singularity}; see, e.g., \cite{kuijlaars2011universality}.
The local statistics at singular points exhibit non-standard universality classes due to the impact of singular points on the surrounding geometry, which can lead to deviations from the typical behaviour observed at regular points, cf.\ \cite{MR4543814}.
Additionally, the insertion of a point charge, also known as the Christoffel perturbation has a~physical application for instance, in the context of the massive quantum field theory \cite{MR2302902,akemann2021skew}.
Let us also mention that the insertion of a point charge has been extensively studied in the context of planar (skew-)orthogonal polynomials, see, e.g., \cite{akemann2021skew,berezin2022planar,lee2020strong} and references therein.
On the one hand, the local statistics of the ensemble also depends on the local geometry of the droplet.
Typically, the focus is on whether the droplet (at the point we zoom in) locally resembles the complex plane or the strip, see, e.g., \cite{MR4030288} and references therein.
The strip regime arises when the particles vary randomly within a thin band of height proportional to the their typical spacing.
In our present case, these regimes can be made by considering the cases where the width of the droplet $S$ in \eqref{droplet} is of order $O(1)$ or $O(1/N)$.
The former is called \emph{strong non-unitarity} and the latter is called \emph{weak non-unitarity} (or \emph{almost-circular regime}).
The latter regime is of particular interest as it generates interpolations between typical one- and two-dimensional statistics.

In summary, we should distinguish the following three different regimes.
\begin{itemize}\itemsep=0pt
 \item[(a)] \textit{At regular points in the limit of strong non-unitarity.} (Cf.\ Figure~\ref{Fig_eigenvalues}\,(a).)
 This means the case where the width of the droplet $S$ in \eqref{droplet} is of order $O(1)$, and the zooming point $p \in \R$ we look at the local statistics is away from the origin.
 To investigate this regime, we set the parameters as
 \begin{equation*}
 L=aN, \qquad n=(b+1)N, \qquad \mbox{with fixed} \quad a,b \ge 0,
 \end{equation*}
 which in the Coulomb gas picture of the Appendix~\ref{A1} corresponds to the external charges at the poles being proportional to $N$.
 Note that with this choice of the parameters, the inner and outer radii in \eqref{droplet} satisfy
 \begin{equation}\label{eq1.13}
 r_1= \sqrt{ \frac{a}{b+1} } +O\left(\frac{1}{N}\right), \qquad r_2 = \sqrt{ \frac{a+1}{ b } }+O\left(\frac{1}{N}\right) \qquad \mbox{as } N \to \infty.
 \end{equation}
 \item[(b)] \textit{At regular points in the limit of weak non-unitarity.} (Cf.\ Figure~\ref{Fig_eigenvalues}\,(b).)
 This means the case where the droplet $S$ is close to the unit circle and its width is of order~$O(1/N)$.
 For this regime, we set
 \begin{equation} \label{L n almost-circular}
 L= \frac{N^2}{\rho^2}-N, \qquad n= \frac{N^2}{\rho^2}, \qquad \mbox{with fixed} \quad \rho>0.
 \end{equation}
 This choice of parameters implies that we impose strong charges (proportional to $N^2$) both at the origin and the infinity which makes the droplet close to the unit circle.
 Indeed, one can see that
\begin{equation*}
	r_1 = 1-\frac{\rho^2}{2N}+O\left(\frac{1}{N^{2}}\right), \qquad r_2 = 1+\frac{\rho^2}{2N}+O\left(\frac{1}{N^{2}}\right), \qquad \mbox{as } N \to \infty.
\end{equation*}
 \item[(c)] \textit{At the singular origin.} (Cf.\ Figure~\ref{Fig_eigenvalues}\,(c).)
 This covers the case where the droplet contains the origin, i.e., $r_1=o(1).$ For this, we set
 \begin{equation*} 
 L>0 \quad \mbox{fixed}, \qquad n=(b+1)N.
 \end{equation*}
 Here the charge at the north pole in the Coulomb gas picture of Appendix~\ref{A1} is $O(1)$.
 Then we have
 \begin{equation*}
 r_1 =O\left(\frac{1}{N}\right) , \qquad r_2 = \frac{1}{\sqrt{b}}+O\left(\frac{1}{N}\right), \qquad \mbox{as } N \to \infty.
 \end{equation*}
\end{itemize}

It is convenient to introduce and recall some notations to describe the scaling limits.
Let us define
\begin{equation} \label{fz}
f_z(u):=\frac12 \erfc\big(\sqrt{2}(z-u)\big).
\end{equation}
Recall that the two-parametric Mittag-Leffler function $E_{a,b}(z)$ is given by
\begin{equation}\label{ML ftn}
E_{a,b}(z):=\sum_{k=0}^\infty \frac{z^k}{\Gamma(ak+b)}.
\end{equation}
We also write $W(f,g):=fg'-gf'$ for the Wronskian.
For a given $p \in \R$, we set
\begin{equation} \label{delta}
\delta := \frac{n+L}{N} \frac{1}{\big(1+p^2\big)^2},
\end{equation}
which corresponds to the density \eqref{density} of the ensemble at the point $p$.
Now we are ready to state our main results. Without loss of generality, it suffices to consider the case $p \ge 0.$

\begin{Theorem}[scaling limits of the eigenvalue correlations] \label{Thm_main}
For a fixed $p \ge 0$, let
\begin{equation} \label{RNk rescaling}
R_{N,k}(z_1,\dots, z_k) := \bfR_{N,k}\left(p+\frac{z_{1}}{\sqrt{N\delta}},\dots,p+\frac{z_{k}}{\sqrt{N\delta}}\right),
\end{equation}
where $\bfR_{N,k}$ and $\delta$ are given by \eqref{bfRNk def} and \eqref{delta}.
Then the following holds.
\begin{itemize}\itemsep=0pt
\item[$(a)$] At regular points in the limit of strong non-unitarity.
Let $L=aN$, $n=(b+1)N$ with fixed $a \ge 0$, $b>0$. Let $p>0$ be fixed. Then as $N\to\infty$,
\begin{align*}
R_{N,k}(z_1,\dots, z_k)&=
 \Pf \Bigg[ {\rm e}^{-|z_j|^2-|z_l|^2}
		\begin{pmatrix}
			\kappa_{(\textup{s})}(z_j,z_l) & \kappa_{(\textup{s})}(z_j,\bar{z}_l)
			 \\
\kappa_{(\textup{s})}(\bar{z}_j,z_l) & \kappa_{(\textup{s})}(\bar{z}_j,\bar{z}_l)
		\end{pmatrix}
		\Bigg]_{j,l=1}^k\\
&\quad \times \prod_{j=1}^k (\bar{z}_j-z_j) + o(1),
\end{align*}
uniformly for $z_{1},\ldots,z_{k}$ in compact subsets of $\mathbb{C}$, where
\begin{gather}
\kappa_{(\textup{s})}(z,w) :=\sqrt{\pi} {\rm e}^{z^2+w^2} \int_{E} W(f_{w},f_{z})(u) \,{\rm d}u,
\nonumber\\
 E= \begin{cases}
(-\infty,\infty) &\textup{if } r_1<p<r_2,
\\
(-\infty,0) &\textup{if } p=r_1 \textup{ or } p=r_2. \label{kappa standard}
\end{cases}
\end{gather}
Here $f_z$ is given by \eqref{fz}.
\item[$(b)$] At regular points in the limit of weak non-unitarity.
Let $L$ and $n$ be given by~\eqref{L n almost-circular}. Let $p=1$. Then as $N\to\infty$,
\begin{align*}
R_{N,k}(z_1,\dots, z_k)
&= \Pf \Bigg[ {\rm e}^{-|z_j|^2-|z_l|^2}
		\begin{pmatrix}
			\kappa_{(\textup{w})}(z_j,z_l) & \kappa_{(\textup{w})}(z_j,\bar{z}_l)
			\\
\kappa_{(\textup{w})}(\bar{z}_j,z_l) & \kappa_{(\textup{w})}(\bar{z}_j,\bar{z}_l)
		\end{pmatrix}
		\Bigg]_{j,l=1}^k\\
&\quad \times\prod_{j=1}^k (\bar{z}_j-z_j) + o(1),
\end{align*}
uniformly for $z_{1},\ldots,z_{k}$ in compact subsets of $\mathbb{C}$, where
\begin{gather}
\kappa_{(\textup{w})}(z,w) :=\sqrt{\pi} {\rm e}^{z^2+w^2} \left( \int_{-a}^a W(f_{w},f_{z})(u) \,{\rm d}u +f_w(a)f_z(-a) -f_z(a)f_w(-a) \right),
\nonumber\\
a= \frac{\rho}{2\sqrt{2}}.\label{kappa almost-circular}
\end{gather}
\item[$(c)$] At the singular origin.
Let $L \ge 0$ be fixed and $n=(b+1)N$ with fixed $b>0$. Let $p=0$. Then as $N\to\infty$,
\begin{align*}
{\rm d} R_{N,k}(z_1,\dots, z_k)
& = \Pf \Bigg[ {\rm e}^{-|z_j|^2-|z_l|^2}
		\begin{pmatrix}
			\kappa_{(\textup{o})}(z_j,z_l) & \kappa_{(\textup{o})}(z_j,\bar{z}_l)
\\
\kappa_{(\textup{o})}(\bar{z}_j,z_l) & \kappa_{(\textup{o})}(\bar{z}_j,\bar{z}_l)
		\end{pmatrix}
		\Bigg]_{j,l=1}^k\\
&\quad \times \prod_{j=1}^k (\bar{z}_j-z_j) + o(1),
\end{align*}
uniformly for $z_{1},\ldots,z_{k}$ in compact subsets of $\mathbb{C}$, where
	\begin{equation} \label{kappa insertion}
		\kappa_{(\textup{o})}(z,w)= 2 (2zw)^{2L} \int_0^1 s^{2L} \big(z{\rm e}^{(1-s^2)z^2}-w{\rm e}^{(1-s^2)w^2}\big) E_{2,1+2L}\big((2szw)^2\big)\,{\rm d}s. 	
	\end{equation}
	Here, $E_{a,b}$ is the two-parametric Mittag-Leffler function \eqref{ML ftn}.
\end{itemize}
\end{Theorem}

The limiting kernel of the form \eqref{kappa standard} was introduced in \cite[Theorem 2.1]{akemann2021scaling} as a scaling limit of the planar symplectic Ginibre ensemble.
Here $E=(-\infty,\infty)$ corresponds to the bulk case, whereas $E=(-\infty,0)$ corresponds to the edge case.
(We also refer to \cite[Remark 2.4]{akemann2021scaling} for more discussions on the role of the integral domain $E$.)
Therefore Theorem~\ref{Thm_main}\,(a) shows that in the limit of strong non-unitarity, the spherical induced symplectic Ginibre ensemble is contained in the universality class of the planar Ginibre ensemble.
Note also that for the bulk case when~$E=(-\infty,\infty)$, the integral in \eqref{kappa standard} can be further simplified, which gives rise to the expression
\begin{equation} \label{kappa standard bulk}
	\kappa_{(\textup{s})}(z,w) = \sqrt{\pi}\,{\rm e}^{z^2+w^2} \erf(z-w) \qquad \mbox{if} \quad E=(-\infty,\infty).
\end{equation}
This form of the kernel appeared in \cite{MR1928853, Mehta}.

The limiting kernel \eqref{kappa almost-circular} was introduced very recently in \cite[Theorem 1.1\,(b)]{byun2022almost} as a scaling limit of the planar induced symplectic ensemble in the almost-circular regime.
Thus Theorem~\ref{Thm_main}\,(b) also establishes the universality in this regime.
An interesting feature of the kernel~\eqref{kappa almost-circular} is that it interpolates the bulk scaling limits of the symplectic Ginibre ensemble which form Pfaffian point processes ($\rho \to \infty$) and those of the chiral Gaussian unitary ensemble ($\rho \to 0$) which form determinantal point processes.
We refer to \cite[Remark 1.4 and Proposition 1.5]{byun2022almost} for more details about this interpolating property.

Finally, the limiting kernel of the form \eqref{kappa insertion} appeared in \cite[Theorem 2.4 and Example 2.6]{akemann2021scaling} (with $c=2L$, $\lambda=1$) as a scaling limit of the planar induced symplectic Ginibre ensemble at the origin having spectral singularity.
Therefore Theorem~\ref{Thm_main}\,(c) again shows the universality and also asserts that under the insertion of a point charge, it is the strength of the charge (i.e., $4L$ in \eqref{Q potential}) that determines the universality class.
We also mention that if $L=0$, one can see from $E_{2,1}\big(z^2\big)= \cosh(z)$ that the kernel \eqref{kappa insertion} agrees with the kernel \eqref{kappa standard bulk}.
Furthermore, it follows from the relation
\begin{align} \label{ML P relation}
		2 E_{2,1+c}\big(z^2\big) = E_{1,1+c}(z)+E_{1,1+c}(-z)
		={\rm e}^z z^{-c} P(c,z)+{\rm e}^{-z} (-z)^{-c} P(c,-z),
	\end{align}
	where
	\begin{equation} \label{incomplete Gamma P}
		P(c,z):=\frac{1}{\Gamma(c)}\int_{0}^{z} t^{c-1} {\rm e}^{-t} \,{\rm d}t, \qquad c>0,
	\end{equation}
	is the (regularised) incomplete gamma function, that we have an alternative representation
	\begin{equation*}
	\begin{split}
	 \kappa_{(\textup{o})}(z,w)=& \int_0^1 \big(z{\rm e}^{(1-s^2)z^2}-w{\rm e}^{(1-s^2)w^2}\big)
	 \\
	 &\quad \times \big( {\rm e}^{2szw} P(2L,2szw)+(-1)^{-2L} {\rm e}^{-2szw} P(2L,-2szw) \big)\,{\rm d}s.
	\end{split}
	\end{equation*}

In Theorem~\ref{Thm_main}, we have focused on the scaling limits along the real axis, i.e., $p \in \R.$
In general, it can be expected that away from the real axis (i.e., $p \in \C \setminus \R$), the scaling limits of the ensemble \eqref{Gibbs} become determinantal with the correlation kernel of its complex counterpart; see~\cite{akemann2019universal} for a heuristic discussion for this.
(Such a statement was shown in \cite{byun2022almost} for the planar induced symplectic ensemble.)
For the spherical induced symplectic Ginibre ensemble, the scaling limits away from the real axis was studied in \cite[Section 6]{MR3612266}, where the authors derived the universal $1$-point functions.

Further points for investigation are also suggested. One is the study of the
so called hole probability, i.e., the probability that a prescribed region is
free of eigenvalues. In the case of the complex Ginibre ensemble, this was
first investigated long ago in \cite{GHS1988}, and in a generalised form
has been the subject of a number of recent works \cite{BC22, Ch21b, Ch21a,fenzl2022precise, lacroix2019intermediate, lacroix2019rotating}.
Another is the study of fluctuation formulas associated with linear statistics;
see the recent review \cite[Section 3.5]{Fo22c} and references therein, and Appendix~\ref{B1} for results
relating to (\ref{Gibbs}) in the case of radial symmetry.

The rest of this paper is organised as follows.
Section~\ref{Section_outline} begins with the finite $N$ result of Proposition~\ref{Prop_finite N} and then identifies a rescaling of the correlation functions valid to leading
order in $N$. Next, in Proposition~\ref{Prop_asymptotic CDI}, the large $N$ form of the differential equation of
Proposition~\ref{Prop_finite N} in the various regimes of interest for Theorem~\ref{Thm_main} is given. The proof of this result is deferred until Section~\ref{Section_remaining proofs}.
The final new result of Section~\ref{Section_outline},
Lemma~\ref{Lem_sol of DEs}, is to present the solutions of the limiting
differential equations. Section~\ref{Section_outline} concludes by showing how these various results can
be assembled to prove Theorem~\ref{Thm_main}. The main content of Section~\ref{Section_remaining proofs} is the proofs of
Propositions~\ref{Prop_finite N} and \ref{Prop_asymptotic CDI}, stated but not proved in earlier sections.

\section{Proof of Theorem~\ref{Thm_main}} \label{Section_outline}

This section culminates in the proof of Theorem~\ref{Thm_main}.
 For reader's convenience, we first present a summary
 of the strategy.
\begin{itemize}\itemsep=0pt
 \item[(i)] In Lemma~\ref{Lem_Rnk basic rescaling}, we first obtain the structure of the correlation function which contains the rescaled skew-kernel $\wt{\kappa}_N$ in \eqref{kappa rescaling} and Gaussian terms.
 This follows from the explicit formula given in Proposition~\ref{Prop_finite N}\,(a).
 \item[(ii)] In Proposition~\ref{Prop_asymptotic CDI}, we derive the asymptotic behaviour of $\pa_z \wt{\kappa}_N(z,w)$.
 For this, we use~\eqref{CDI} in Proposition~\ref{Prop_finite N}\,(b) and compute the asymptotic expansions of its inhomogeneous terms.
 \item[(iii)] In Lemma~\ref{Lem_sol of DEs}, we solve the differential equation appearing in Proposition~\ref{Prop_asymptotic CDI}, which gives rise to the explicit formulas of the limiting correlation kernels in Theorem~\ref{Thm_main}.
\end{itemize}

As already mentioned in the final paragraph of the above section, the proofs of Propositions~\ref{Prop_finite N} and \ref{Prop_asymptotic CDI} are given separately in Section~\ref{Section_remaining proofs}.
We begin with deriving the basic structure of the correlation functions using Propo\-si\-tion~\ref{Prop_finite N}\,(a).

\begin{Lemma}[structure of the correlation functions] \label{Lem_Rnk basic rescaling}
For a fixed $p \ge 0$, let
\begin{equation} \label{kappa rescaling}
\wt{\kappa}_N(z,w):= \frac{1}{\big(1+p^2\big)^3\,(N\delta)^{3/2}} \wt{\bfkappa}_N\left(p+\frac{z}{\sqrt{N\delta}},p+\frac{w}{\sqrt{N}}\right),
\end{equation}
where $\wt{\bfkappa}_N$ is given by \eqref{bfkappa tilda G wh} and $\delta>0$ is given by \eqref{delta}.
Recall that $R_{N,k}$ is given by \eqref{RNk rescaling}.
Then under the same assumptions on Theorem~{\rm \ref{Thm_main}} and in each case $(a)$, $(b)$, $(c)$, we have
\begin{align*}
R_{N,k}(z_1,\dots, z_k)&
 = \Pf \Bigg[
{\rm e}^{-|z_j|^2-|z_l|^2}
	\begin{pmatrix}
		{\rm e}^{z_j^2+z_l^2} \wt{\kappa}_N(z_j,z_l) &
	{\rm e}^{z_j^2+\bar{z}_l^2} \wt{\kappa}_N(z_j,\bar{z}_l)
		\\
	{\rm e}^{\bar{z}_j^2+z_l^2}		\wt{\kappa}_N(\bar{z}_j,z_l) & 	{\rm e}^{\bar{z}_j^2+\bar{z}_l^2} \wt{\kappa}_N(\bar{z}_j,\bar{z}_l)
	\end{pmatrix} \Bigg]_{ j,l=1 }^k\\
&\quad\times \prod_{j=1}^{k} (\bar{z}_j-z_j) + O\left(\frac{1}{\sqrt{N}}\right),
\end{align*}
uniformly for $z_{1},\ldots,z_{k}$ in compact subsets of $\mathbb{C}$, as $N \to \infty$.
\end{Lemma}

\begin{proof}
Recall that the weight function $\omega$ is given by \eqref{omega}.
Then it follows from direct computations that in each case (a), (b), (c) of Theorem~\ref{Thm_main},
\begin{equation} \label{omega asymp}
 \omega\left( p+\frac{z}{ \sqrt{N \delta} } \right) = \big(1+p^2\big)^{-\frac32} {\rm e}^{ -|z|^2+\frac12(z^2+\bar{z}^2) } + O\left(\frac{1}{\sqrt{N}}\right), \qquad \mbox{as } N\to\infty.
\end{equation}
More precisely,
by \eqref{omega}, we have
\begin{align*}
 \omega\left( p+\frac{z}{\sqrt{N\delta}} \right) &= \big(1+p^2\big)^{-\frac32}
 \left| 1+ \frac{2p}{1+p^2} \frac{z}{ \sqrt{N\delta} }+ \frac{1}{1+p^2} \frac{z^2}{N\delta} \right|^{ n+L-\frac12 } \\
 &\quad \times \left( 1+ \frac{p}{1+p^2} \frac{z+\bar{z}}{ \sqrt{N \delta } } + \frac{1}{1+p^2} \frac{|z|^2}{ N\delta } \right)^{-n-L-1}.
\end{align*}
Furthermore, in each case (a), (b), (c), we have
\begin{align*}
\log \left| 1+ \frac{2p}{1+p^2} \frac{z}{ \sqrt{N\delta} }+ \frac{1}{1+p^2} \frac{z^2}{N\delta} \right| & = \frac{p}{1+p^2} \frac{z+\bar{z}}{ \sqrt{N\delta} } + \frac{\big(1-p^2\big)\big(z^2+\bar{z}^2\big) }{ 2\big(1+p^2\big)^2 } \frac{ 1 }{ N\delta } +O\big(N^{-3/2}\big)
\\
&= \frac{p(z+\bar{z})}{ \sqrt{n+L} } + \frac{\big(1-p^2\big)\big(z^2+\bar{z}^2\big)}{ 2 } \frac{ 1 }{ n+L } +O\big(N^{-3/2}\big)
\end{align*}
and
\begin{align*}
\log \left(\! 1+ \frac{p}{1+p^2} \frac{z+\bar{z}}{ \sqrt{N \delta } } + \frac{1}{1+p^2} \frac{|z|^2}{ N\delta } \right)& = \frac{p}{1+p^2} \frac{z+\bar{z}}{ \sqrt{N\delta} } + \frac{2|z|^2\!-p^2\big(z^2+\bar{z}^2\big) }{ 2\big(1+p^2\big)^2 } \frac{ 1 }{ N\delta } +O\big(N^{-3/2}\big)
\\
&= \frac{p(z+\bar{z})}{ \sqrt{n+L} } + \frac{2|z|^2\!-p^2\big(z^2+\bar{z}^2\big) }{ 2 } \frac{ 1 }{ n+L } +O\big(N^{-3/2}\big),
\end{align*}
as $N \to \infty,$ where we have used \eqref{delta}.
Combining the above, we obtain \eqref{omega asymp}.

For given $p \ge 0$, we write
\begin{equation} \label{zeta z}
\zeta_j= p + \frac{z_j}{\sqrt{N\delta}},
\end{equation}
where $\delta$ is given by \eqref{delta}.
Then by combining Proposition~\ref{Prop_finite N} with \eqref{RNk rescaling}, \eqref{omega asymp} and \eqref{zeta z}, we obtain
\begin{align*}
R_{N,k}(z_1,\dots, z_k)&= \frac{1}{(N\delta)^k} \bfR_{N,k}(\zeta_1,\dots, \zeta_k)
\\
&= \Pf \Bigg[ \frac{
{\rm e}^{-|z_j|^2-|z_l|^2} }{ \big(1+p^2\big)^{3} (N\delta)^{3/2} }
	\begin{pmatrix}
		{\rm e}^{z_j^2+z_l^2} \wt{\bfkappa}_N(\zeta_j,\zeta_l) &
	{\rm e}^{z_j^2+\bar{z}_l^2} \wt{\bfkappa}_N(\zeta_j,\bar{\zeta}_l)
			\\
	{\rm e}^{\bar{z}_j^2+z_l^2}		\wt{\bfkappa}_N(\bar{\zeta}_j,\zeta_l) & 	{\rm e}^{\bar{z}_j^2+\bar{z}_l^2} \wt{\bfkappa}_N(\bar{\zeta}_j,\bar{\zeta}_l)
	\end{pmatrix} \Bigg]_{ j,l=1 }^k \\
&\quad\times \prod_{j=1}^{k} (\bar{z}_j-z_j) + O\left(\frac{1}{\sqrt{N}}\right).
\end{align*}
Here, we have used the fact that the Pfaffian of a correlation kernel is invariant under the multiplication by cocycles, see, e.g., \cite[p.~19]{akemann2021scaling}.
Lemma~\ref{Lem_Rnk basic rescaling} follows from \eqref{kappa rescaling}.
\end{proof}

The next step is to derive the asymptotic behaviour of the derivative of $\wt{\kappa}_N$ in \eqref{kappa rescaling}.
This step crucially relies on Proposition~\ref{Prop_finite N}\,(b).

\begin{Proposition}[large-$N$ expansions of the differential equations] \label{Prop_asymptotic CDI}
As $N\to \infty$, the following hold.
\begin{itemize}\itemsep=0pt
 \item[$(a)$] Under the setup of Theorem~{\rm \ref{Thm_main}}\,$(a)$,
 \begin{equation} \label{ODE a}
 \pa_z \wt{\kappa}_N(z,w)= F_{\textup{(s)}}(z,w)+o(1),
 \end{equation}
 uniformly for $z$, $w$ in compact subsets of $\mathbb{C}$, where
 \begin{equation} \label{Fs inhomo}
F_{\textup{(s)}}(z,w):=
\begin{cases}
\displaystyle 2\, {\rm e}^{-(z-w)^2} & \text{if } r_1<p<r_2,
\\
\displaystyle {\rm e}^{-(z-w)^2} \erfc(z+w) -\frac{ {\rm e}^{-2z^2} }{ \sqrt{2} } \erfc\big(\sqrt{2}w\big) & \text{if } p=r_1,r_2.
\end{cases}
 \end{equation}
 \item[$(b)$] Under the setup of Theorem~{\rm \ref{Thm_main}}\,$(b)$,
 \begin{equation} \label{ODE b}
 \pa_z \wt{\kappa}_N(z,w)= F_{\textup{(w)}}(z,w)+o(1),
 \end{equation}
 uniformly for $z$, $w$ in compact subsets of $\mathbb{C}$, where
 \begin{align}
& F_{\textup{(w)}}(z,w) := {\rm e}^{-(z-w)^2} \big( \erfc\big(z+w-\tfrac{\rho}{\sqrt{2}}\big)- \erfc\big(z+w+\tfrac{\rho}{\sqrt{2}}\big) \big)
\nonumber\\
& \quad - \frac{1}{\sqrt{2}} \big( {\rm e}^{ -(\sqrt{2}z-\frac{\rho}{2})^2 }+ {\rm e}^{ -(\sqrt{2}z+\frac{\rho}{2})^2 } \big)\big( \erfc\big(\sqrt{2}w-\tfrac{\rho}{2}\big)-\erfc\big(\sqrt{2}w+\tfrac{\rho}{2}\big) \big).\label{Fw inhomo}
 \end{align}
\item[$(c)$] Under the setup of Theorem~{\rm \ref{Thm_main}}\,$(c)$,
 \begin{equation} \label{ODE c}
 \pa_z \wt{\kappa}_N(z,w)= F_{\textup{(o)}}(z,w)+o(1),
 \end{equation}
 uniformly for $z$, $w$ in compact subsets of $\mathbb{C}$, where
 \begin{align} \label{Fo inhomo}
 F_{\textup{(o)}}(z,w) :=2 {\rm e}^{-(z-w)^2} P(2L,2zw)- \frac{2\sqrt{\pi}}{ \Gamma(L) } z^{2L-1} {\rm e}^{-z^2} P\big(L+\tfrac12, w^2\big) .
\end{align}
Here $P$ is the regularised incomplete gamma function \eqref{incomplete Gamma P}.
\end{itemize}
\end{Proposition}

Let us mention that the case $L=0$ in \eqref{Fo inhomo} can be interpreted by using the fact $1/\Gamma(k+1) = 0$ for a negative integer $k.$
The proof of this proposition will be given in the next section.
Finally, we solve the differential equations appearing in Proposition~\ref{Prop_asymptotic CDI}.
The following lemma is an immediate consequence of several results established in \cite{akemann2021scaling, byun2022almost}.

\begin{Lemma} \label{Lem_sol of DEs}
Let
\begin{gather}
\mathcal{K}_{\textup{(s)}}(z,w) := {\rm e}^{-z^2-w^2}\kappa_{\textup{(s)}}(z,w), \label{mathcal K kappa s}
\\
\mathcal{K}_{\textup{(w)}}(z,w) := {\rm e}^{-z^2-w^2}\kappa_{\textup{(w)}}(z,w),
\label{mathcal K kappa w}
\\
\mathcal{K}_{\textup{(o)}}(z,w) := {\rm e}^{-z^2-w^2}\kappa_{\textup{(o)}}(z,w),
\label{mathcal K kappa o}
\end{gather}
where $\kappa_{\textup{(s)}}$, $\kappa_{\textup{(w)}}$ and $\kappa_{\textup{(o)}}$ are given by \eqref{kappa standard}, \eqref{kappa almost-circular} and \eqref{kappa insertion}. Then the following hold.
\begin{itemize}\itemsep=0pt
 \item[$(a)$] For a given $w \in \C$, the function $z \mapsto \mathcal{K}_{\textup{(s)}}(z,w)$ is a unique solution to
 \begin{equation*}
 \pa_z \mathcal{K}_{\textup{(s)}}(z,w)= F_{\textup{(s)}}(z,w), \qquad \mathcal{K}_{\textup{(s)}}(z,w)|_{z=w}=0,
 \end{equation*}
 where $F_{\textup{(s)}}$ is given by \eqref{Fs inhomo}.
 \item[$(b)$] For a given $w \in \C$, the function $z \mapsto \mathcal{K}_{\textup{(w)}}(z,w)$ is a unique solution to
 \begin{equation*}
 \pa_z \mathcal{K}_{\textup{(w)}}(z,w)= F_{\textup{(w)}}(z,w), \qquad \mathcal{K}_{\textup{(w)}}(z,w)|_{z=w}=0,
 \end{equation*}
 where $F_{\textup{(w)}}$ is given by \eqref{Fw inhomo}.

 \item[$(c)$] For a given $w \in \C$, the function $z \mapsto \mathcal{K}_{\textup{(o)}}(z,w)$ is a unique solution to
 \begin{equation*}
 \pa_z \mathcal{K}_{\textup{(o)}}(z,w)= F_{\textup{(o)}}(z,w),\qquad \mathcal{K}_{\textup{(o)}}(z,w)|_{z=w}=0,
 \end{equation*}
 where $F_{\textup{(o)}}$ is given by \eqref{Fo inhomo}.
\end{itemize}
\end{Lemma}

\begin{proof}By \eqref{kappa standard bulk}, the first assertion (a) for the bulk case when $E=(-\infty,\infty)$ is trivial.
For the edge case when $E=(-\infty,0)$, this was shown in \cite[p.~21]{akemann2021scaling}.

The second assertion (b) was shown in \cite[Section~3.2]{byun2022almost}.

Finally, the third assertion (c) follows from \cite[Section 4]{akemann2021scaling}.
More precisely, the equation in the statement (c) is a special case of \cite[equation~(4.4)]{akemann2021scaling} with $\lambda=1, c=2L$ up to a trivial transformation. Here, we also use the relation~\eqref{ML P relation}.
Then this differential equation was solved in \cite[Section~4.2]{akemann2021scaling}.
\end{proof}

Let us now combine the results introduced above and finish the proof of Theorem~\ref{Thm_main}.

\begin{proof}[Proof of Theorem~\ref{Thm_main}]
For a given $w$, we view \eqref{ODE a}, \eqref{ODE b}, \eqref{ODE c} as first-order ordinary differential equations in $z$ with initial conditions $\wt{\kappa}_N(z,w)|_{z=w}=0.$
Combining Proposition~\ref{Prop_asymptotic CDI}, Lemma~\ref{Lem_sol of DEs} and \cite[Lemma 3.10]{byun2021universal}, we obtain that
\begin{equation} \label{wt kappa asymp}
\wt{\kappa}_N(z,w)= {\rm e}^{-z^2-w^2} \begin{cases}
\kappa_{ \textup{s} }(z,w)+o(1) &\textup{for the case (a)},
\\
\kappa_{ \textup{w} }(z,w)+o(1) &\textup{for the case (b)},
\\
\kappa_{ \textup{o} }(z,w)+o(1) &\textup{for the case (c)}.
\end{cases}
\end{equation}
Here we also have used \eqref{mathcal K kappa s}, \eqref{mathcal K kappa w} and \eqref{mathcal K kappa o}.
Furthermore, note that both $\wt{\kappa}_N$ and the $o(1)$-terms in \eqref{wt kappa asymp} are anti-symmetric in $z$ and $w$.
In particular, the entire proof remains valid if the roles of $z$ and $w$ are interchanged.
The theorem now follows from Lemma~\ref{Lem_Rnk basic rescaling} and \eqref{wt kappa asymp}.
\end{proof}

\section{Proof of Propositions~\ref{Prop_finite N} and~\ref{Prop_asymptotic CDI}} \label{Section_remaining proofs}

In this section we present the proofs of Propositions~\ref{Prop_finite N} and~\ref{Prop_asymptotic CDI}. Both these results have been used in the proof of Theorem~\ref{Thm_main}.

\subsection{Skew-orthogonal polynomial kernels}

Using the general theory on planar symplectic ensembles and skew-orthogonal polynomials, we first show the following lemma.

Recall that the potential $Q$ and the correlation function $\bfR_{N,k}$ are given by \eqref{Q potential} and \eqref{bfRNk def}.

\begin{Lemma}
\label{Lem_Rnk basic}
We have
\begin{gather} \label{bfR Pfa basic}
	\bfR_{N,k}(\zeta_1,\dots, \zeta_k) = \Pf \Bigg[
	{\rm e}^{ -N(Q(\zeta_j)+Q(\zeta_l)) }
	\begin{pmatrix}
		\bfkappa_N(\zeta_j,\zeta_l) & \bfkappa_N(\zeta_j,\bar{\zeta}_l)
		\\
		\bfkappa_N(\bar{\zeta}_j,\zeta_l) & \bfkappa_N(\bar{\zeta}_j,\bar{\zeta}_l)
	\end{pmatrix} \Bigg]_{ j,l=1 }^k \prod_{j=1}^{k} \big(\bar{\zeta}_j-\zeta_j\big),
\end{gather}
where
\begin{equation} \label{bfkappaN GN}
	\begin{split}
		\bfkappa_N(\zeta,\eta):=\boldsymbol{G}_N(\zeta,\eta)-\boldsymbol{G}_N(\eta,\zeta),
	\end{split}
\end{equation}
and
\begin{equation*}
	\begin{split} 
		\boldsymbol{G}_N(\zeta,\eta)&:=
		 \pi \frac{ \Gamma(2n+2L+2) }{ 2^{2L+2n+1} }
		 \\
		 & \quad \times \sum_{k=0}^{N-1} \sum_{l=0}^{k} \frac{ \zeta^{2k+1} \eta^{2l} }{ \Gamma\big(k+L+\frac32\big) \Gamma(n-k) \Gamma\big(n-l+\frac12\big) \Gamma(l+L+1) } .
	\end{split}
\end{equation*}
\end{Lemma}

We stress that Lemma~\ref{Lem_Rnk basic} is an immediate consequence of \cite[Proposition 4]{MR3066630}.
Nevertheless, as this lemma is crucially used in the present work, we briefly recall the proof.

\begin{proof}[Proof of Lemma~\ref{Lem_Rnk basic}]
First, let us consider the ensemble \eqref{Gibbs} with a general potential $Q$.
Define the skew-symmetric form
\begin{equation*}
\langle f, g \rangle_s := \int_{\C} \big( f(\zeta) g\big(\bar{\zeta}\big) - g(\zeta) f\big(\bar{\zeta}\big) \big) \big(\zeta - \bar{\zeta}\big) {\rm e}^{-2N Q(\zeta)} \,{\rm d}A(\zeta).
\end{equation*}
Then the skew-orthogonal polynomial $q_m$ of degree $m$ is defined by the condition: for all $k, l \in \mathbb{N}$
\begin{equation*}
\langle q_{2k}, q_{2l} \rangle_s = \langle q_{2k+1}, q_{2l+1} \rangle_s = 0, \qquad \langle q_{2k}, q_{2l+1} \rangle_s = -\langle q_{2l+1}, q_{2k} \rangle_s = r_k \,\delta_{k, l}.
\end{equation*}
Here, $\delta_{k, l}$ is the Kronecker delta.
Then it is well known \cite{MR1928853} that the correlation function \eqref{bfRNk def} is of the form \eqref{bfR Pfa basic} with the canonical skew-kernel
\begin{gather}
\bfkappa_N(\zeta,\eta):=\sum_{k=0}^{N-1} \frac{q_{2k+1}(\zeta) q_{2k}(\eta) -q_{2k}(\zeta) q_{2k+1}(\eta)}{r_k},
\nonumber\\
\boldsymbol{G}_N(\zeta,\eta):= \sum_{k=0}^{N-1} \frac{q_{2k+1}(\zeta) q_{2k}(\eta) }{r_k}.\label{bfkappaN skewOP}
\end{gather}
Thus it suffices to compute the skew-orthogonal polynomials.

Let us now consider a general radially symmetric potential $Q(\zeta)=Q(|\zeta|)$. We write
\begin{equation*}
 h_k:=\int_\C |\zeta|^{2k} {\rm e}^{-2N Q(\zeta)} \,{\rm d}A(\zeta)
\end{equation*}
for the (squared) orthogonal norm.
Then it follows from \cite[Corollary~3.2]{akemann2021skew} that
\begin{gather} \label{skew op_rad}
	q_{2k+1}(\zeta)=\zeta^{2k+1}, \qquad
	q_{2k}(\zeta)=\zeta^{2k}+\sum_{l=0}^{k-1} \zeta^{2l} \prod_{j=0}^{k-l-1} \frac{h_{2l+2j+2} }{ h_{2l+2j+1} }, \qquad
	r_k=2h_{2k+1}.
\end{gather}

We now turn to the potential \eqref{Q potential}. For this
\begin{equation*}
h_k = 2 \int_0^\infty \frac{ r^{2k+4L+1} }{ (1+r^2)^{2(n+L+1)} } \,{\rm d}r = \frac{ \Gamma(k+2L+1) \Gamma(2n-k+1) }{ \Gamma(2n+2L+2) },
\end{equation*}
and so by using \eqref{skew op_rad}, we obtain
\begin{gather}
q_{2k}(\zeta)= \frac{ \Gamma(k+L+1) } { \Gamma\big(k-n+\frac12\big) } \sum_{l=0}^{k} (-1)^{k-l} \frac{ \Gamma\big(l-n+\frac12\big) }{ \Gamma(l+L+1) } \zeta^{2l},
\nonumber\\
 r_k= \frac{ 2 \Gamma(2k+2L+2) \Gamma(2n-2k) }{ \Gamma(2n+2L+2) }.\label{sop spherical}
\end{gather}
(Note that these formulas were also derived \cite[Proposition~4]{MR3066630}.)
Combining \eqref{bfkappaN skewOP}, \eqref{sop spherical} and the basic functional relations
\begin{equation*} 
\Gamma(z)\Gamma(1-z)=\frac{\pi}{\sin(\pi z)},\qquad \Gamma(2z)=\frac{2^{2z-1}}{\sqrt{\pi}} \Gamma(z)\Gamma\big(z+\tfrac12\big)
\end{equation*}
of the gamma function, we obtain
\begin{align*}
\boldsymbol{G}_N(\zeta,\eta) &=\sum_{k=0}^{N-1} \frac{ \Gamma(2n+2L+2) }{ 2 \Gamma(2k+2L+2) \Gamma(2n-2k) } \frac{ \Gamma(k+L+1) } { \Gamma\big(k-n+\frac12\big) } \\
&\quad \times \sum_{l=0}^{k} (-1)^{k-l} \frac{ \Gamma\big(l-n+\frac12\big) }{ \Gamma(l+L+1) } \zeta^{2k+1} \eta^{2l}
 \\
 &=\sum_{k=0}^{N-1} \frac{ \Gamma(2n+2L+2) \Gamma\big(n-k+\frac12\big) }{ 2 \Gamma(2k+2L+2) \Gamma(2n-2k) } \frac{ \Gamma(k+L+1) } { \pi (-1)^{k-n} } \\
 &\quad \times \sum_{l=0}^{k} (-1)^{k-l} \frac{ \Gamma\big(l-n+\frac12\big) }{ \Gamma(l+L+1) } \zeta^{2k+1} \eta^{2l}
 \\
 &=\Gamma(2n+2L+2) \sum_{k=0}^{N-1} \frac{ 1 }{ 2^{2L+2n+1} \Gamma\big(k+L+\frac32\big) \Gamma(n-k) }\\
 &\quad \times \sum_{l=0}^{k} (-1)^{n-l} \frac{ \Gamma\big(l-n+\frac12\big) }{ \Gamma(l+L+1) } \zeta^{2k+1} \eta^{2l}.
\end{align*}
This completes the proof.
\end{proof}

\subsection{Proof of Propositions~\ref{Prop_finite N}}

Let
\begin{equation} \label{bfkappaN GN hat}
		\wh{\bfkappa}_N(\zeta,\eta):=\wh{\boldsymbol{G}}_N(\zeta,\eta)-\wh{\boldsymbol{G}}_N(\eta,\zeta) =\big( \big(1+\zeta^2\big) \big(1+\eta^2\big) \big)^{n+L-\frac12} \wt{\bfkappa}_N(\zeta,\eta),
\end{equation}
where $\wh{\boldsymbol{G}}_N$ and $\wt{\bfkappa}_N$ are given by \eqref{GN hat} and \eqref{bfkappa tilda G wh}.

The key step to prove Proposition~\ref{Prop_finite N}\,(b) is the following lemma.

\begin{Lemma}\label{Lem_CDI pre}
We have
\begin{align*}
 \big(1+\zeta^2\big)\pa_\zeta \wh{\bfkappa}_N(\zeta,\eta) & = 2 \zeta \left( n+L-\frac12 \right) \wh{\bfkappa}_N(\zeta,\eta)
 \\
 &\quad + \frac{ \Gamma(2n+2L+2) }{2} \sum_{k=0}^{2N-1} \frac{ \zeta^{k+2L} \eta^{k+2L} }{ \Gamma(k+2L+1) \Gamma(2n-k) }
\\
& \quad - \frac{ 	 \pi \, \Gamma(2n+2L+2) \, \zeta^{2N+2L} }{ 2^{2L+2n} \Gamma\big(N+L+\frac12\big) \Gamma(n-N) } \sum_{k=0}^{N-1} \frac{ \eta^{2k+2L} }{ \Gamma\big(n-k+\frac12\big) \Gamma(k+L+1) }
\\
& \quad - \frac{ 	 \pi \, \Gamma(2n+2L+2) }{ 2^{2L+2n} \Gamma\big(n+\frac12\big) \Gamma(L) } \zeta^{2L-1} \sum_{k=0}^{N-1} \frac{ \eta^{2k+2L+1} }{ \Gamma(n-k) \Gamma\big(k+L+\frac32\big) } .
\end{align*}
\end{Lemma}

\begin{proof}Let us first compute $ \pa_\zeta \wh{\boldsymbol{G}}_N(\zeta,\eta)$.
Note that
\begin{align*}
&\pa_\zeta \sum_{k=0}^{N-1} \sum_{l=0}^{k} \frac{ \zeta^{2k+2L+1} \eta^{2l+2L} }{ \Gamma\big(k+L+\frac32\big) \Gamma(n-k) \Gamma\big(n-l+\frac12\big) \Gamma(l+L+1) }
\\
&\qquad= 2 \zeta \sum_{k=0}^{N-1} \sum_{l=0}^{k} \frac{ \zeta^{2k+2L-1} \eta^{2l+2L} }{ \Gamma\big(k+L+\frac12\big) \Gamma(n-k) \Gamma\big(n-l+\frac12\big) \Gamma(l+L+1) }
\\
&\qquad= 2 \frac{ \zeta^{2L} \eta^{2L} }{ \Gamma\big(L+\frac12\big) \Gamma(n) \Gamma\big(n+\frac12\big) \Gamma(L+1) }
\\
&\qquad\quad + 2 \zeta \sum_{k=0}^{N-2} \sum_{l=0}^{k+1} \frac{ \zeta^{2k+2L+1} \eta^{2l+2L} }{ \Gamma\big(k+L+\frac32\big) \Gamma(n-k-1) \Gamma\big(n-l+\frac12\big) \Gamma(l+L+1) }.
\end{align*}
Here, we have
\begin{align*}
&2 \zeta\sum_{k=0}^{N-2} \sum_{l=0}^{k+1} \frac{ \zeta^{2k+2L+1} \eta^{2l+2L} }{ \Gamma\big(k+L+\frac32\big) \Gamma(n-k-1) \Gamma\big(n-l+\frac12\big) \Gamma(l+L+1) }
\\
&\qquad= 2 \zeta\sum_{k=0}^{N-1} \sum_{l=0}^{k} \frac{ \zeta^{2k+2L+1} \eta^{2l+2L} }{ \Gamma\big(k+L+\frac32\big) \Gamma(n-k-1) \Gamma\big(n-l+\frac12\big) \Gamma(l+L+1) }
\\
&\qquad\quad - 2 \zeta \sum_{l=0}^{N-1} \frac{ \zeta^{2N+2L-1} \eta^{2l+2L} }{ \Gamma\big(N+L+\frac12\big) \Gamma(n-N) \Gamma\big(n-l+\frac12\big) \Gamma(l+L+1) }
\\
&\qquad\quad + 2 \sum_{k=1}^{N-1} \frac{ \zeta^{2k+2L} \eta^{2k+2L} }{ \Gamma\big(k+L+\frac12\big) \Gamma(n-k) \Gamma\big(n-k+\frac12\big) \Gamma(k+L+1) }.
\end{align*}
Therefore we obtain
\begin{align*}
&\pa_\zeta \sum_{k=0}^{N-1} \sum_{l=0}^{k} \frac{ \zeta^{2k+2L+1} \eta^{2l+2L} }{ \Gamma\big(k+L+\frac32\big) \Gamma(n-k) \Gamma\big(n-l+\frac12\big) \Gamma(l+L+1) }
\\
&\qquad= 2 \zeta\sum_{k=0}^{N-1} \sum_{l=0}^{k} \frac{ \zeta^{2k+2L+1} \eta^{2l+2L} }{ \Gamma\big(k+L+\frac32\big) \Gamma(n-k-1) \Gamma\big(n-l+\frac12\big) \Gamma(l+L+1) }
\\
&\qquad\quad - 2 \sum_{l=0}^{N-1} \frac{ \zeta^{2N+2L} \eta^{2l+2L} }{ \Gamma\big(N+L+\frac12\big) \Gamma(n-N) \Gamma\big(n-l+\frac12\big) \Gamma(l+L+1) }
\\
&\qquad\quad + 2 \sum_{k=0}^{N-1} \frac{ \zeta^{2k+2L} \eta^{2k+2L} }{ \Gamma\big(k+L+\frac12\big) \Gamma(n-k) \Gamma\big(n-k+\frac12\big) \Gamma(k+L+1) }.
\end{align*}
Note here that by \eqref{GN hat}, we have
\begin{align*}
&\pi \frac{ \Gamma(2n+2L+2) }{ 2^{2L+2n+1} } \sum_{k=0}^{N-1} \sum_{l=0}^{k} \frac{ \zeta^{2k+2L+1} \eta^{2l+2L} }{ \Gamma\big(k+L+\frac32\big) \Gamma(n-k-1) \Gamma\big(n-l+\frac12\big) \Gamma(l+L+1) }
\\
&\qquad= 	 \pi \frac{ \Gamma(2n+2L+2) }{ 2^{2L+2n+1} } \sum_{k=0}^{N-1} \sum_{l=0}^{k} \frac{ (n-k-1) \zeta^{2k+2L+1} \eta^{2l+2L} }{ \Gamma\big(k+L+\frac32\big) \Gamma(n-k) \Gamma\big(n-l+\frac12\big) \Gamma(l+L+1) }
\\
&\qquad= 	 \pi \frac{ \Gamma(2n+2L+2) }{ 2^{2L+2n+1} } \sum_{k=0}^{N-1} \sum_{l=0}^{k} \frac{ \big( \big(n+L-\frac12\big)-\big(k+L+\frac12\big) \big) \zeta^{2k+2L+1} \eta^{2l+2L} }{ \Gamma\big(k+L+\frac32\big) \Gamma(n-k) \Gamma\big(n-l+\frac12\big) \Gamma(l+L+1) }
\\
&\qquad= \left(n+L-\frac12\right) \wh{\boldsymbol{G}}_N(\zeta,\eta) -\frac{1}{2} \zeta \pa_\zeta \wh{\boldsymbol{G}}_N(\zeta,\eta).
\end{align*}
Combining all of the above equations, we conclude
\begin{align}
\label{GN deri 1}
 \pa_\zeta \wh{\boldsymbol{G}}_N(\zeta,\eta) &= 2 \zeta \left(n+L-\frac12\right) \wh{\boldsymbol{G}}_N(\zeta,\eta) - \zeta^2 \pa_\zeta \wh{\boldsymbol{G}}_N(\zeta,\eta)
\\
&\quad -	 \pi \frac{ \Gamma(2n+2L+2) }{ 2^{2L+2n} } \sum_{l=0}^{N-1} \frac{ \zeta^{2N+2L} \eta^{2l+2L} }{ \Gamma\big(N+L+\frac12\big) \Gamma(n-N) \Gamma\big(n-l+\frac12\big) \Gamma(l+L+1) }
\nonumber\\
&\quad +	 \pi \frac{ \Gamma(2n+2L+2) }{ 2^{2L+2n} } \sum_{k=0}^{N-1} \frac{ \zeta^{2k+2L} \eta^{2k+2L} }{ \Gamma\big(k+L+\frac12\big) \Gamma(n-k) \Gamma\big(n-k+\frac12\big) \Gamma(k+L+1) }.\nonumber
\end{align}

Next, we compute $ \pa_\zeta \wh{\boldsymbol{G}}_N(\eta,\zeta)$.
By similar computations as above, we have
\begin{align*}
&\pa_\zeta \sum_{k=0}^{N-1} \sum_{l=0}^{k} \frac{ \eta^{2k+2L+1} \zeta^{2l+2L} }{ \Gamma\big(k+L+\frac32\big) \Gamma(n-k) \Gamma\big(n-l+\frac12\big) \Gamma(l+L+1) }
\\
&\qquad=2 \sum_{k=0}^{N-1} \frac{ \eta^{2k+2L+1} \zeta^{2L-1} }{ \Gamma\big(k+L+\frac32\big) \Gamma(n-k) \Gamma\big(n+\frac12\big) \Gamma(L) }
\\
&\qquad\quad + 2 \zeta \sum_{k=0}^{N-1} \sum_{l=0}^{k-1} \frac{ \eta^{2k+2L+1} \zeta^{2l+2L} }{ \Gamma\big(k+L+\frac32\big) \Gamma(n-k) \Gamma\big(n-l-\frac12\big) \Gamma(l+L+1) } .
\end{align*}
Here, the last term is rearranged as
\begin{align*}
& \sum_{k=0}^{N-1} \sum_{l=0}^{k-1} \frac{ \eta^{2k+2L+1} \zeta^{2l+2L} }{ \Gamma\big(k+L+\frac32\big) \Gamma(n-k) \Gamma\big(n-l-\frac12\big) \Gamma(l+L+1) }
\\
&\qquad= \sum_{k=0}^{N-1} \sum_{l=0}^{k} \frac{ \eta^{2k+2L+1} \zeta^{2l+2L} }{ \Gamma\big(k+L+\frac32\big) \Gamma(n-k) \Gamma\big(n-l-\frac12\big) \Gamma(l+L+1) }
\\
&\qquad\quad - \sum_{k=0}^{N-1} \frac{ \eta^{2k+2L+1} \zeta^{2k+2L} }{ \Gamma\big(k+L+\frac32\big) \Gamma(n-k) \Gamma\big(n-k-\frac12\big) \Gamma(k+L+1) } .
\end{align*}
This gives rise to
\begin{align*}
&\pa_\zeta \sum_{k=0}^{N-1} \sum_{l=0}^{k} \frac{ \eta^{2k+2L+1} \zeta^{2l+2L} }{ \Gamma\big(k+L+\frac32\big) \Gamma(n-k) \Gamma\big(n-l+\frac12\big) \Gamma(l+L+1) }
\\
&\qquad= 2 \zeta \sum_{k=0}^{N-1} \sum_{l=0}^{k} \frac{ \eta^{2k+2L+1} \zeta^{2l+2L} }{ \Gamma\big(k+L+\frac32\big) \Gamma(n-k) \Gamma\big(n-l-\frac12\big) \Gamma(l+L+1) }
\\
&\qquad \quad + 2 \sum_{k=0}^{N-1} \frac{ \eta^{2k+2L+1} \zeta^{2L-1} }{ \Gamma\big(k+L+\frac32\big) \Gamma(n-k) \Gamma\big(n+\frac12\big) \Gamma(L) }
\\
&\qquad \quad - 2 \zeta \sum_{k=0}^{N-1} \frac{ \eta^{2k+2L+1} \zeta^{2k+2L} }{ \Gamma\big(k+L+\frac32\big) \Gamma(n-k) \Gamma\big(n-k-\frac12\big) \Gamma(k+L+1) } .
\end{align*}
By using \eqref{GN hat}, we also have
\begin{align*}
&\pi \frac{ \Gamma(2n+2L+2) }{ 2^{2L+2n+1} } \sum_{k=0}^{N-1} \sum_{l=0}^{k} \frac{ \eta^{2k+2L+1} \zeta^{2l+2L} }{ \Gamma\big(k+L+\frac32\big) \Gamma(n-k) \Gamma\big(n-l-\frac12\big) \Gamma(l+L+1) }
\\
&\qquad = 	 \pi \frac{ \Gamma(2n+2L+2) }{ 2^{2L+2n+1} } \sum_{k=0}^{N-1} \sum_{l=0}^{k} \frac{ \big(n-l-\frac12\big) \eta^{2k+2L+1} \zeta^{2l+2L} }{ \Gamma\big(k+L+\frac32\big) \Gamma(n-k) \Gamma\big(n-l+\frac12\big) \Gamma(l+L+1) }
\\
& \qquad = 	 \pi \frac{ \Gamma(2n+2L+2) }{ 2^{2L+2n+1} } \sum_{k=0}^{N-1} \sum_{l=0}^{k} \frac{ \big(\big(n+L-\frac12\big)-(l+L)\big) \eta^{2k+2L+1} \zeta^{2l+2L} }{ \Gamma\big(k+L+\frac32\big) \Gamma(n-k) \Gamma\big(n-l+\frac12\big) \Gamma(l+L+1) }
\\
&\qquad = \left( n+L-\frac12 \right) \wh{\boldsymbol{G}}_N(\eta,\zeta) -\frac12 \zeta \pa_\zeta \wh{\boldsymbol{G}}_N(\eta,\zeta).
\end{align*}
Therefore we have shown that
\begin{align}
 \pa_\zeta \wh{\boldsymbol{G}}_N(\eta,\zeta) &= 2 \zeta \left( n+L-\frac12 \right) \wh{\boldsymbol{G}}_N(\eta,\zeta) - \zeta^2 \pa_\zeta \wh{\boldsymbol{G}}_N(\eta,\zeta)\label{GN deri 2}
 \\
&\quad +	 \pi \frac{ \Gamma(2n+2L+2) }{ 2^{2L+2n} } \sum_{k=0}^{N-1} \frac{ \eta^{2k+2L+1} \zeta^{2L-1} }{ \Gamma\big(k+L+\frac32\big) \Gamma(n-k) \Gamma\big(n+\frac12\big) \Gamma(L) }\nonumber
\\
&\quad -	 \pi \frac{ \Gamma(2n+2L+2) }{ 2^{2L+2n} } \sum_{k=0}^{N-1} \frac{ \eta^{2k+2L+1} \zeta^{2k+2L+1} }{ \Gamma\big(k+L+\frac32\big) \Gamma(n-k) \Gamma\big(n-k-\frac12\big) \Gamma(k+L+1) } .\nonumber
\end{align}
The lemma follows from \eqref{GN deri 1}, \eqref{GN deri 2} and \eqref{bfkappaN GN hat}.
\end{proof}

We now finish the proof of Proposition~\ref{Prop_finite N}.

\begin{proof}[Proof of Proposition~\ref{Prop_finite N}]
We first show the first part.
By combining Lemma~\ref{Lem_Rnk basic} with~\eqref{Q potential}, \eqref{bfkappa tilda G wh} and \eqref{bfkappaN GN}, we obtain
\begin{align*}
& \bfR_{N,k}(\zeta_1,\dots, \zeta_k) = \prod_{j=1}^{k} (\bar{\zeta}_j-\zeta_j)\, \Pf \Bigg[ \frac{1}{ \big( 	\big(1+|\zeta_j|^2\big) \big(1+|\zeta_l|^2\big) \big)^{n+L+1} }
	\\
	& 	\times\begin{pmatrix}
	\big(\big(1+\zeta_j^2\big)\big(1+\zeta_l^2\big)\big)^{n+L-\frac12} 	\wt{\bfkappa}_N(\zeta_j,\zeta_l) &
	\big(\big(1+\zeta_j^2\big)\big(1+\bar{\zeta}_l^2\big)\big)^{n+L-\frac12} \wt{\bfkappa}_N(\zeta_j,\bar{\zeta}_l)
		\\
\big(\big(1+\bar{\zeta}_j^2\big)\big(1+\zeta_l^2\big)\big)^{n+L-\frac12} 	\wt{\bfkappa}_N(\bar{\zeta}_j,\zeta_l) & \big(\big(1+\bar{\zeta}_j^2\big)\big(1+\bar{\zeta}_l^2\big)\big)^{n+L-\frac12} \wt{\bfkappa}_N(\bar{\zeta}_j,\bar{\zeta}_l)
	\end{pmatrix} \Bigg]_{ j,l=1 }^k .
\end{align*}
Then \eqref{bfRNk bfkappa tilde} follows from the basic properties of Pfaffians.
For instance, we have
\begin{equation*}
\bfR_{N,1}(\zeta)= \frac{ \big|1+\zeta^2\big|^{2n+2L-1} }{ \big(1+|\zeta|^2\big)^{2n+2L+2} } \wt{\bfkappa}_N\big(\zeta,\bar{\zeta}\big)\big(\bar{\zeta}-\zeta\big)
 \end{equation*}
and{\samepage
\begin{align*}
\begin{split}
\bfR_{N,2}(\zeta,\eta)& = \frac{ \big|1+\zeta^2\big|^{2n+2L-1} \big|1+\eta^2\big|^{2n+2L-1} }{ \big(1+|\zeta|^2\big)^{2n+2L+2} \big(1+|\eta|^2\big)^{2n+2L+2} }
\\
&\quad \times \big( \wt{\bfkappa}_N\big(\zeta,\bar{\zeta}\big) \wt{\bfkappa}_N(\eta,\bar{\eta})- |\wt{\bfkappa}_N(\zeta,\eta)|^2 + |\wt{\bfkappa}_N(\zeta,\bar{\eta})|^2 \big) \big(\bar{\zeta}-\zeta\big) (\bar{\eta}-\eta).
\end{split}
\end{align*}
This establishes Proposition~\ref{Prop_finite N}\,(a).}

For the second assertion, recall that $\mathfrak{p}$ and $\mathfrak{q}$ are given in \eqref{p q mathfrak} and that
\begin{equation} \label{binomial gamma}
\binom{r}{k}= \frac{\Gamma(r+1)}{ \Gamma(k+1) \Gamma(r-k+1) }.
\end{equation}
Then after some straightforward computations using Lemma~\ref{Lem_CDI pre}, \eqref{binomial gamma} and the transform \eqref{bfkappaN GN hat}, the desired formula \eqref{CDI} follows.
\end{proof}

\subsection{Proof of Proposition~\ref{Prop_asymptotic CDI}}

It remains to show Proposition~\ref{Prop_asymptotic CDI} to validate the proof of Theorem~\ref{Thm_main}.
We begin with the following lemma which is the rescaled version of Proposition~\ref{Prop_finite N}\,(b).

\begin{Lemma}\label{Lem_CDI rescaling}We have
\begin{equation}\label{CDI rescaled}
\pa_z \wt{\kappa}_N(z,w) = \RN{1}_N^{(1)}(z,w)\RN{1}_N^{(2)}(z,w)-\RN{2}_N^{(1)}(z,w)\RN{2}_N^{(2)}(z,w)-\RN{3}_N^{(1)}(z,w)\RN{3}_N^{(2)}(z,w),
\end{equation}
where
\begin{align}
&\RN{1}_N^{(1)}(z,w) := \frac{1}{\big(1+p^2\big)^3(N\delta)^2} \frac{1}{1+\zeta^2} \frac{\big(1+\zeta \eta\big)^{2n+2L-1} }{\big(\big(1+\zeta^2\big)\big(1+\eta^2\big)\big)^{n+L-\frac12} } \frac{ \Gamma(2n+2L+2) }{2 \Gamma(2n+2L)} , \label{RN1 wt 1}
\\
&\RN{2}_N^{(1)}(z,w) := \frac{1}{\big(1+p^2\big)^3(N\delta)^2}\frac{ \zeta^{2N+2L}}{ \big(1+\zeta^2\big)^{n+L+\frac12} } \frac{ 	 \pi \, \Gamma(2n+2L+2)/2^{2L+2n} }{ \Gamma\big(N+L+\frac12\big) \Gamma(n-N) \Gamma\big(n+L+\frac12\big) } , \label{RN2 wt 1}
\\
&\RN{3}_N^{(1)}(z,w) := \frac{1}{\big(1+p^2\big)^3}\frac{1}{(N\delta)^2} \frac{\zeta^{2L-1}}{ \big(1+\zeta^2\big)^{n+L+\frac12} } \frac{ 	 \pi \, \Gamma(2n+2L+2)/2^{2L+2n} }{ \Gamma\big(n+\frac12\big) \Gamma(L) \Gamma\big(n+L+\frac12\big) }, \label{RN3 wt 1}
\end{align}
and
\begin{align}
& \RN{1}_N^{(2)}(z,w) :=\sum_{k=0}^{2N-1} \binom{2n+2L-1}{k+2L} \left( \frac{\zeta \eta}{1+\zeta \eta}\right)^{k+2L} \left( \frac{1}{ 1+\zeta \eta } \right)^{2n-k-1}, \label{RN1 wt 2}
\\
& \RN{2}_N^{(2)}(z,w) :=\sum_{k=0}^{N-1} \binom{n+L-\frac12}{k+L} \left( \frac{\eta^2}{1+\eta^2}\right)^{k+L} \left( \frac{1}{1+\eta^2}\right)^{n-k-\frac12}, \label{RN2 wt 2}
\\
& \RN{3}_N^{(2)}(z,w) := \sum_{k=0}^{N-1} \binom{n+L-\frac12}{k+L+\frac12} \left( \frac{\eta^2}{1+\eta^2}\right)^{k+L+\frac12} \left( \frac{1}{1+\eta^2} \right)^{ n-k-1 }. \label{RN3 wt 2}
\end{align}
Here,
\begin{equation} \label{zeta eta z w}
\zeta= p+\frac{z}{ \sqrt{N \delta} } , \qquad \eta= p+\frac{w}{ \sqrt{N \delta} }.
\end{equation}
\end{Lemma}

\begin{proof}This is an immediate consequence of \eqref{kappa rescaling} and \eqref{CDI}.
\end{proof}

We need to analyse the right-hand side of \eqref{CDI rescaled}.
For this, we need the following lemma.

\begin{Lemma}\label{Lem_CDI asymp 1}
Recall that $\RN{1}_N^{(1)}$, $\RN{2}_N^{(1)}$, $\RN{3}_N^{(1)}$ are given by \eqref{RN1 wt 1}, \eqref{RN2 wt 1}, \eqref{RN3 wt 1} and $\zeta, \eta$ are given by~\eqref{zeta eta z w}.
Let $\epsilon>0$ be a small constant.
As $N\to \infty$, the following hold.
\begin{itemize}\itemsep=0pt
 \item[$(a)$] Under the setup of Theorem~{\rm \ref{Thm_main}}\,$(a)$, we have
\begin{align*}
&\RN{1}_N^{(1)}(z,w) = 2 {\rm e}^{-(z-w)^2} +O\left(\frac{1}{\sqrt{N}}\right),
\\
&\RN{2}_N^{(1)}(z,w) =
\begin{cases}
\sqrt{2}\, {\rm e}^{-2z^2} &\text{if } p=r_2,
\\
O\big({\rm e}^{-N \epsilon}\big) &\text{otherwise},
\end{cases}
\\
&\RN{3}_N^{(1)}(z,w) =
\begin{cases}
\sqrt{2}\, {\rm e}^{-2z^2} &\text{if } p=r_1,
\\
O\big({\rm e}^{-N \epsilon}\big) &\text{otherwise},
\end{cases}
\end{align*}
uniformly for $z$, $w$ in compact subsets of $\mathbb{C}$.

 \item[$(b)$] Under the setup of Theorem~{\rm \ref{Thm_main}}\,$(b)$, we have
\begin{align*}
& \RN{1}_N^{(1)}(z,w) = 2\, {\rm e}^{-(z-w)^2} +O\left(\frac{1}{N}\right),
\\
& \RN{2}_N^{(1)}(z,w) = \sqrt{2} \, {\rm e}^{- (\sqrt{2}z-\frac{\rho}{2} )^2} +O\left(\frac{1}{N}\right),
\\
& \RN{3}_N^{(1)}(z,w) = \sqrt{2} \, {\rm e}^{- (\sqrt{2}z+\frac{\rho}{2} )^2} +O\left(\frac{1}{N}\right),
\end{align*}
 uniformly for $z$, $w$ in compact subsets of $\mathbb{C}$.

 \item[$(c)$] Under the setup of Theorem~{\rm \ref{Thm_main}}\,$(c)$, we have
\begin{align*}
& \RN{1}_N^{(1)}(z,w) = 2\, {\rm e}^{-(z-w)^2} +O\left(\frac{1}{N}\right),
\\
& \RN{2}_N^{(1)}(z,w) = O\big( {\rm e}^{-N \epsilon} \big),
\\
& \RN{3}_N^{(1)}(z,w) = \frac{2\sqrt{\pi}}{ \Gamma(L) } z^{2L-1} {\rm e}^{-z^2} +O\left(\frac{1}{N}\right),
\end{align*}
 uniformly for $z$, $w$ in compact subsets of $\mathbb{C}$.
\end{itemize}
 \end{Lemma}

\begin{proof}This follows from long but straightforward computations repeatably using Stirling's formula.
The most non-trivial part is the computations for $\RN{2}_N^{(1)}$ and $\RN{3}_N^{(1)}$ under the setup of Theorem~\ref{Thm_main}\,(a).
In this case, we have
\begin{align*}
\RN{2}_N^{(1)}(z,w) &= \left( \frac{ (1+a+b)^{1+a+b} }{ (1+a)^{1+a} b^b }\frac{ p^{2+2a} }{ \big(1+p^2\big)^{1+a+b} } \right)^N \exp \left( \frac{ 2\big(1+a-bp^2\big) }{ p \sqrt{1+a+b} }\,z \sqrt{N} \right)
\\
&\quad \times \Bigg[ \exp \left( - \frac{1+a+(3+3a+b)p^2-bp^4 }{ (1+a+b)p^2 } z^2 \right) \frac{\big(1+p^2\big)^{1/2}}{(1+a+b)^{1/2}} \sqrt{2b} +O\left(\frac{1}{\sqrt{N}}\right)\Bigg]
\end{align*}
and
\begin{align*}
\RN{3}_N^{(1)}(z,w) & = \left( \frac{ (1+a+b)^{a+b+1} }{ a^{a} (1+b)^{1+b} }\frac{ p^{2a} }{ \big(1+p^2\big)^{1+a+b} } \right)^N \exp \left( \frac{ 2\big(a-(1+b)p^2\big) }{ p \sqrt{1+a+b} }\,z \sqrt{N} \right)
\\
&\quad \times \Bigg[ \exp\!\left(\! -\frac{ a+(1+3a+b)p^2-(1+b)p^4 }{ (1+a+b)p^2 } z^2 \!\right)\! \frac{\big(1+p^2\big)^{1/2}}{(1+a+b)^{1/2}} \frac{ \sqrt{2a} }{ p} +O\!\left(\!\frac{1}{\sqrt{N}}\!\right)\!\Bigg]
\end{align*}
as $N\to\infty$. Then the desired asymptotic expansion follows from the asymptotic formulas for~$r_1$ and~$r_2$ given in~\eqref{eq1.13}.
\end{proof}

\begin{Lemma}\label{Lem_CDI asymp 2}
Recall that $\RN{1}_N^{(2)}$, $\RN{2}_N^{(2)}$, $\RN{3}_N^{(2)}$ are given by \eqref{RN1 wt 2}, \eqref{RN2 wt 2}, \eqref{RN3 wt 2} and $\zeta, \eta$ are given by~\eqref{zeta eta z w}.
As $N\to \infty$, the following holds.
\begin{itemize}\itemsep=0pt
 \item[$(a)$] Under the setup of Theorem~{\rm \ref{Thm_main}}\,$(a)$, we have
\begin{align}
& \RN{1}_N^{(2)}(z,w) = \begin{cases}
1+o(1) &\text{if } r_1<p<r_2
\\
\frac12\erfc(z+w) +o(1) &\text{if } p=r_1,r_2,
\end{cases} \label{RN 1 2 asymp a}
\\
& \RN{2}_N^{(2)}(z,w) = \begin{cases}
1+o(1) &\text{if } r_1<p<r_2,
\\
\frac12\erfc\big(\sqrt{2}w\big)+o(1) &\text{if } p=r_1,r_2 ,
\end{cases} \label{RN 2 2 asymp a}
\\
& \RN{3}_N^{(2)}(z,w) = \begin{cases}
1+o(1) &\text{if } r_1<p<r_2,
\\
\frac12\erfc\big(\sqrt{2}w\big)+o(1) &\text{if } p=r_1,r_2 ,
\end{cases} \label{RN 3 2 asymp a}
\end{align}
uniformly for $z$, $w$ in compact subsets of $\mathbb{C}$.
\item[$(b)$] Under the setup of Theorem~{\rm \ref{Thm_main}}\,$(b)$, we have
\begin{align}
& \RN{1}_N^{(2)}(z,w) = \frac12 \bigl( \erfc\big(z+w-\tfrac{\rho}{\sqrt{2}}\big)- \erfc\big(z+w+\tfrac{\rho}{\sqrt{2}}\bigr) \big)+o(1), \label{RN 1 2 asymp b}
\\
& \RN{2}_N^{(2)}(z,w) = \frac12 \bigl( \erfc\big(\sqrt{2}w-\tfrac{\rho}{2}\big)-\erfc\big(\sqrt{2}w+\tfrac{\rho}{2}\big) \bigr)+o(1) , \label{RN 2 2 asymp b}
\\
& \RN{3}_N^{(2)}(z,w) = \frac12 \bigl( \erfc\big(\sqrt{2}w-\tfrac{\rho}{2}\big)-\erfc\big(\sqrt{2}w+\tfrac{\rho}{2}\big) \bigr)+o(1), \label{RN 3 2 asymp b}
\end{align}
uniformly for $z$, $w$ in compact subsets of $\mathbb{C}$.

\item[$(c)$] Under the setup of Theorem~{\rm \ref{Thm_main}}\,$(c)$, we have
\begin{align}
& \RN{1}_N^{(2)}(z,w) = P(2L,2zw)+o(1), \label{RN 1 2 asymp c}
\\
& \RN{2}_N^{(2)}(z,w) = P\big(L,w^2\big)+o(1) , \label{RN 2 2 asymp c}
\\
& \RN{3}_N^{(2)}(z,w) = P \big(L+\tfrac12,w^2\big)+o(1), \label{RN 3 2 asymp c}
\end{align}
uniformly for $z$, $w$ in compact subsets of $\mathbb{C}$. Here $P$ is the regularised incomplete gamma function.
\end{itemize}
 \end{Lemma}

\begin{proof}Recall that $\mathfrak{p}$ and $\mathfrak{q}$ are given by \eqref{p q mathfrak}.
We first present a probabilistic proof of the lemma, which requires in particular that $z,w \in \mathbb R$.
This is instructive as it clearly shows the appearance of the $\erfc$ and the incomplete gamma functions in the context of the normal and Poisson approximations of the binomial distributions. We then extend the validity to complex~$z$,~$w$ by Vitali's theorem.
The setting of the latter is a sequence of uniformly bounded analytic functions of a single complex variable in a region $\Omega$. Vitali's theorem gives that the convergence of the sequence proved on a dense subset of $\Omega$ can be extended to convergence on all compact subsets of $\Omega$.
Such a strategy has been applied in related settings in, e.g.,~\cite[proof of Theorem 3.3]{BG99}.

Let $X \sim B(2n+2L-1,\mathfrak{p})$ be the binomial distribution assuming that $2n+2L-1$, $2L$ are integers and $z,w \in \R$.
Then $\RN{1}_N^{(2)}$ can be rewritten as
\begin{gather*}
\RN{1}_N^{(2)}(z,w) = \Prob ( 2L \le X \le 2N+2L-1 )
\\
\hphantom{\RN{1}_N^{(2)}(z,w)}= \Prob\Bigg( \frac{2L-(2n+2L-1)\mathfrak{p}}{ \sqrt{ (2n+2L-1)\mathfrak{p}(1-\mathfrak{p}) } } \le \frac{X-(2n+2L-1)\mathfrak{p}}{ \sqrt{ (2n+2L-1)\mathfrak{p}(1-\mathfrak{p}) } }
\\
\hphantom{\RN{1}_N^{(2)}(z,w)= \Prob\Bigg( \frac{2L-(2n+2L-1)\mathfrak{p}}{ \sqrt{ (2n+2L-1)\mathfrak{p}(1-\mathfrak{p}) } }} \le \frac{2N+2L-1 -(2n+2L-1)\mathfrak{p} }{ \sqrt{ (2n+2L-1)\mathfrak{p}(1-\mathfrak{p}) } } \Bigg).
\end{gather*}
Let us consider the setup of Theorem~\ref{Thm_main}\,(a).
Then as $N\to \infty$,
\begin{align*}
\frac{2N+2L-1 -(2n+2L-1)\mathfrak{p} }{ \sqrt{ (2n+2L-1)\mathfrak{p}(1-\mathfrak{p}) } } &= \frac{ \sqrt{2}\big(1+a-bp^2\big) }{ p \sqrt{ 1+a+b } } \sqrt{N}
\\
&\quad- \frac{ \big(1+p^2\big) \big(1+a+bp^2\big) }{ \sqrt{2} (1+a+b) p^2 } (z+w)+O\left(\frac{1}{\sqrt{N}}\right),
\\
\frac{2L-(2n+2L-1)\mathfrak{p}}{ \sqrt{ (2n+2L-1)\mathfrak{p}(1-\mathfrak{p}) } } &= \frac{ \sqrt{2}\big(a-(1+b)p^2\big) }{ p \sqrt{ 1+a+b } } \sqrt{N}
\\
 &\quad- \frac{ \big(1+p^2\big) \big(a+(1+b)p^2\big) }{ \sqrt{2} (1+a+b) p^2 } (z+w)+O\left(\frac{1}{\sqrt{N}}\right),
\end{align*}
uniformly for $z$, $w$ in compact subsets of $\mathbb{C}$.
Recall that $r_1 \sim \sqrt{ \frac{a}{b+1} }$ and $ r_2 \sim \sqrt{ \frac{a+1}{ b } }$.
Then by the Gaussian approximation of the binomial distribution, as $N\to \infty$, we obtain
\begin{equation*}
\Prob ( 2L \le X \le 2N+2L-1 ) \sim
\begin{cases}
\Prob ( -\infty \le Z \le \infty ) &\text{if } r_1<p<r_2,
\\
\Prob ( -\infty \le Z \le -\sqrt{2}(z+w) ) &\text{if } p=r_1,r_2,
\end{cases}
\end{equation*}
where $Z$ is the standard normal distribution.
This gives rise to the desired asymptotic behaviour~\eqref{RN 1 2 asymp a}.
All other asymptotic formulas \eqref{RN 2 2 asymp a}, \eqref{RN 3 2 asymp a}, \eqref{RN 1 2 asymp b}, \eqref{RN 2 2 asymp b}, \eqref{RN 3 2 asymp b} involving the $\erfc$ function follow along the same lines.

Under the setup of Theorem~\ref{Thm_main}\,(c), we have
\begin{equation*}
\mathfrak{p}= \frac{zw}{1+b} \frac{1}{N}+O\left(\frac{1}{N^2}\right), \qquad \mbox{as } N \to \infty.
\end{equation*}
Thus the binomial distribution $X$ is approximated by the Poisson distribution with intensity $\lambda=(2n+2L-1)\mathfrak{p} \sim 2 zw$.
Since the regularised incomplete gamma function is the cumulative distribution function of the Poisson distribution, we have
\begin{equation*}
\Prob ( 2L \le X \le 2N+2L-1 ) \sim P(2L,2zw), \qquad \mbox{as } N \to \infty,
\end{equation*}
which leads to \eqref{RN 1 2 asymp c}.
The other asymptotics \eqref{RN 2 2 asymp c}, \eqref{RN 3 2 asymp c} follow in a similar way.

We now turn to the case with general parameters.
In general, the functions $\RN{1}_N^{(2)}$, $\RN{2}_N^{(2)}$, $\RN{3}_N^{(2)}$ can be written in terms of the (regularised) incomplete beta function
\begin{equation} \label{incomplete beta}
 I_x(a,b):= \frac{ \Gamma(a+b) }{ \Gamma(a) \Gamma(b) } \int_0^x t^{a-1}(1-t)^{b-1}\,{\rm d}t
\end{equation}
as
\begin{align}
& \RN{1}_N^{(2)}(z,w) = I_{\mathfrak{p}}(2L,2n)- I_{\mathfrak{p}}(2N+2L,2n-2N) , \label{RN1 wt 2 beta}
\\
& \RN{2}_N^{(2)}(z,w) = I_{\mathfrak{q}}\big(L,n+\tfrac12\big)- I_{\mathfrak{q}}\big(N+L,n-N+\tfrac12\big), \label{RN2 wt 2 beta}
\\
& \RN{3}_N^{(2)}(z,w) = I_{\mathfrak{q}}\big(L+\tfrac12,n\big)- I_{\mathfrak{q}}\big(N+L+\tfrac12,n-N\big). \label{RN3 wt 2 beta}
\end{align}
For integer valued cases, the expressions \eqref{RN1 wt 2 beta}, \eqref{RN2 wt 2 beta}, \eqref{RN3 wt 2 beta} easily follow from
\begin{equation*}
I_x (m,n-m+1) = \sum_{j=m}^n \binom{n}{j} x^j (1-x)^{n-j},
\end{equation*}
see, e.g., \cite[equation~(8.17.5)]{olver2010nist}.
In general, these follow from the definition of the hypergeometric function in series form
\begin{equation*}
 {}_2 F_1(a,b;c;z)= \frac{\Gamma(c)}{ \Gamma(a)\Gamma(b) } \sum_{s=0}^\infty \frac{ \Gamma(a+s)\Gamma(b+s) }{ \Gamma(c+s) s! } z^s,
\end{equation*}
and the relation
\begin{gather*}
 I_x(a,b) = \frac{ \Gamma(a+b) }{ \Gamma(a) \Gamma(b) } \frac{x^a (1-x)^{b-1}}{ a }\, {}_2 F_1\left(1,1-b ; a+1 ; \frac{x}{x-1} \right).
\end{gather*}

The significance of (\ref{RN1 wt 2 beta})--(\ref{RN3 wt 2 beta}) is that $I_x(a,b)$ as defined by (\ref{incomplete beta}) is an analytic function of $x$ in the cut complex-$x$ plane $\mathbb C \backslash (-\infty, 0)$ provided $|x| < 1$. Appealing to Vitali's theorem on uniform convergence inside of a domain $\mathcal C$ for sequences of analytic functions on $\mathcal C$ (see \cite{Ti64}) then allows the result proved for $z$, $w$ real to be extended to compact sets of the complex plane. The required uniform bound is a consequence of the scaling of these variables by $\sqrt{N}$ as required by (\ref{zeta eta z w}), ensuring that for $z$, $w$ in compact subsets of $\mathbb C$, the limiting sequence remains in the domain of analyticity.
\end{proof}

\begin{Remark}
For $z$, $w$ real, the uniform asymptotic expansions of the incomplete beta function~\eqref{incomplete beta} can be found in \cite[Section~8.18]{olver2010nist} and \cite[Section~11.3.3]{MR1376370}. A method to extend these to the complex plane using a direct argument can be found in \cite[Section~5]{Te79}.
 \end{Remark}

We now finish the proof of Proposition~\ref{Prop_asymptotic CDI}.

\begin{proof}[Proof of Proposition~\ref{Prop_asymptotic CDI}]
This immediately follows by substituting the asymptotic expansions in Lemmas~\ref{Lem_CDI asymp 1} and~\ref{Lem_CDI asymp 2} into the identity \eqref{CDI rescaled}.
\end{proof}

\appendix
\section{Appendix}\label{A1}

Consider an eigenvalue probability density function for $2N$ eigenvalues in the complex plane, coordinates
$\{ \zeta_j \}_{j=1}^N$, specified by
\begin{equation}\label{A.1}
\frac{1}{Z_{2N}} \prod_{1 \le j < k \le 2N} | \zeta_k - \zeta_j |^2 \prod_{l=1}^{2N} {\rm e}^{- N Q(|\zeta_l|)}.
\end{equation}
Here the radial potential~$Q(|\zeta_l|)$ is to be taken as general, subject only to the normalisation~$Z_{2N}$
being well defined. Suppose next that in this functional form only $\{ \zeta_j \}_{j=1}^N$ are independent, with
$\zeta_{j+N} = \bar{\zeta}_j$ ($j=1,\dots,N$). Then (\ref{A.1}) reduces to a probability density function for $N$ eigenvalues
specified by
\begin{equation}\label{A.2}
\frac{1}{\tilde{Z}_{N}} \prod_{1 \le j < k \le N} | \zeta_k - \zeta_j |^2 \prod_{l=1}^{N} {\rm e}^{- 2N Q(|\zeta_l|)}.
\end{equation}
We see that \eqref{matrix model} relates through this prescription to (\ref{A.1}) with $ Q(|\zeta_l|)$ given by \eqref{Q potential}. In this appendix, following \cite{FF11}, we revise the interpretation of (\ref{A.1}) in terms of the Boltzmann factor for a certain one-component two-dimensional Coulomb system, features of which are then inherited by (\ref{A.2}).

The first point to note is the mapping of (\ref{A.1}) from the complex plane to the unit diameter Riemann sphere specified by
\begin{equation}\label{A.2a}
z = {\rm e}^{{\rm i} \phi} \tan(\theta/2), \qquad 0 \le \phi \le 2 \pi, \qquad 0 \le \theta \le \pi,
\end{equation}
which geometrically corresponds to a stereographic projection from the south pole. Introducing the Cayley--Klein parameters
\begin{equation*}
u = \cos(\theta/2) {\rm e}^{{\rm i} \phi/2}, \qquad v = - {\rm i} \sin(\theta/2) {\rm e}^{-{\rm i} \phi/2},
\end{equation*}
and with ${\rm d} S$ denoting an element of the surface area of the sphere, a straightforward calculation shows
\begin{align}
& \prod_{l=1}^m \left( \frac{ |z_l |^2}{1 + | z_l |^2} \right)^{q_1 m} \frac{1}{\big(1 + | z_l |^2\big)^{q_2 m + m + 1} } \prod_{1 \le j < k \le m} | z_k - z_j |^2 \,
{\rm d} z_1 \cdots {\rm d} z_m
\nonumber\\
&\qquad{}= \prod_{l=1}^m | v_l|^{2 q_1 {m} } | u_l|^{2 q_2 {m} } \prod_{1 \le j < k \le m} | u_k v_j - u_j v_k |^2 \,{\rm d} S_1 \cdots {\rm d} S_m,\label{A.3}
\end{align}
where $m := 2 N$. The relevance of (\ref{A.3}) is that with $Q(|\zeta_l|)$ given by \eqref{Q potential}, the left-hand side of~(\ref{A.3}) results with
\begin{equation}\label{A.3a}
q_1 = \frac{2L}{m}, \qquad q_2 = \frac{(2n - m)}{m}.
\end{equation}

The parameters $q_1$, $q_2$ have an electrostatic interpretation on the right-hand side of (\ref{A.3}). This comes about by
first recalling the fact that two points $(\theta, \phi)$ and $(\theta', \phi')$ on a sphere of unit diameter, the solution of the charge
neutral Poisson equation at $(\theta, \phi)$ due to a unit charge at~$(\theta', \phi')$, is (see, e.g.,~\cite[equation~(15.108)]{forrester2010log}) given in terms of
the corresponding Cayley--Klein parameters by
\begin{equation}\label{A.4}
\Phi((\theta,\phi), (\theta', \phi')) = - \log | u' v - u v'|.
\end{equation}
Let there be $m$ unit charges with coordinates $(\theta, \phi)$ interacting pairwise on the sphere via the potential (\ref{A.4}).
This gives an energy
\begin{equation*}
U_0 = - \sum_{1 \le j < k \le m} \log | u_j v_k - u_k v_j|.
\end{equation*}
Suppose that at the north pole there is a fixed charge $m q_1$, and at the south pole there is a fixed charge $m q_2$. The interaction
with the mobile unit charges gives an energy
\begin{equation*}
U_1 = - m q_1 \sum_{j=1}^m \log |v_j| - m q_2 \sum_{j=1}^m \log | u_j|.
\end{equation*}
We see that forming the Boltzmann factor ${\rm e}^{-\beta (U_0 + U_1)}$ gives, for $\beta = 2$, precisely the right-hand side of (\ref{A.3}).

A spherical cap about the north pole with azimuthal angle $\theta$ has surface area $\pi \sin^2 (\theta/2)$. With a charge $m q_1$ at the north
pole, the value of $\theta$, $\theta_{q_1}$ say, which corresponds to a uniform neutralising background charge $-m q_1$ in the spherical cap is such that
\begin{equation*}
\sin^2 (\theta_{q_1}/2) = \frac{q_1}{q_1 + q_2 + 1}.
\end{equation*}
Here the left-hand side is the proportion of the total surface area of the sphere which is in the spherical cap.
On the right-hand side the ratio is obtained by dividing the charge at the north pole by the total charge. Mapped to the
complex plane using (\ref{A.2a}), this gives a radius $r_{q_1}$ such that $r_{q_1}^2 = q_1/(1 + q_2)$ --- note that
this corresponds to $r_1$ in \eqref{droplet}. An analogous calculation for the spherical cap about the charge $m q_2$ at the south pole
corresponding to a charge neutral region leads to $r_{q_2}^2 = (1+q_1)/q_2$, which corresponds to $r_2$ in \eqref{droplet}.

\section{Appendix}\label{B1}

Here consideration is given to fluctuation formulas for linear statistics relating to (\ref{Gibbs}). A linear statistic is the random function
$B := \sum_{j=1}^N b(\zeta_j)$. The mean $\mu_{N,B}$ can be expressed in terms of the eigenvalue density
($1$-point correlation) $\bfR_{N,1}(\zeta)$ according to
\begin{equation*}
\mu_{N,B} = \int_{\mathbb C} b(\zeta) \bfR_{N,1}(\zeta) \,{\rm d} A(\zeta).
\end{equation*}
The variance can be expressed in terms of the one and two point correlations according to
\begin{align}
\sigma_{N,B}^2 & = \int_{\mathbb C} {\rm d} A(\zeta_1) \, b(\zeta_1) \int_{\mathbb C} {\rm d}A(\zeta_2) \, b(\zeta_2)
\nonumber\\
& \quad \times \big ( \bfR_{N,2}(\zeta_1, \zeta_2) - \bfR_{N,1}(\zeta_1) \bfR_{N,1}(\zeta_2) + \bfR_{N,1}(\zeta_1) \delta(\zeta_1 - \zeta_2) \big ).
\label{B.1}
\end{align}

The most appropriate scaling regime to analyse a linear statistic is the macroscopic limit.
The density then has the large $N$ form given by \eqref{density}, supported on the region $S$ specified by~\eqref{droplet}, and so in this setting
\begin{equation}\label{B.2}
\mu_{N,B} \sim (n + L) \int_S \frac{b(\zeta)}{\big(1 + | \zeta|^2\big)^2} \, {\rm d} A(\zeta).
\end{equation}
In particular, this shows the mean is extensive, being proportional to $N$.
In contrast, in the macroscopic limit the variance is expected to be independent of $N$, under the assumption that~$b$ is smooth. The full distribution is expected to be a Gaussian. Heuristic reasoning from the Coulomb gas viewpoint underlying these predictions can be found, e.g., in~\cite[Section~14.4]{forrester2010log}.

The limit formulas for the correlations functions of Theorem~\ref{Thm_main} relate to local rather than global scaling. Upon global scaling the functional form relating to the correlations in (\ref{B.1}) is not expected to be well defined as a
function, but rather to take the form of a distribution; see~\cite[Section 15.4]{forrester2010log}. In fact in the particular case that~$b$ is
smooth and a function of the distance from the origin only, it is possible to compute the limiting form of the variance
indirectly, by considering the large~$N$ form of the characteristic function.

\begin{Proposition}
Consider the radially symmetric linear statistic $B = \sum_{j=1}^N b(|\zeta_j|)$ in relation to the induced symplectic induced spherical ensemble as specified by \eqref{Gibbs} and \eqref{Q potential}. Let $\hat{P}_{N,B}(k)$ denote
the corresponding characteristic function. We have that for large $N$
\begin{equation}\label{B.3}
\hat{P}_{N,B}(k) = {\rm e}^{{\rm i} k \tilde{\mu}_{N,B} - k^2 \tilde{\sigma}_B^2/2 + o(1)},
\end{equation}
 where $ \tilde{\mu}_{N,B}$ denotes the right-hand side of \eqref{B.2} and with $S$ defined as in \eqref{B.2}
\begin{equation}\label{B.4}
 \tilde{\sigma}_B^2 = \frac{1}{8} \int_S || \nabla b(|\zeta|) ||^2 \,{\rm d} A(\zeta).
 \end{equation}
 \end{Proposition}

\begin{proof}[Sketch of proof] For the most part we follow the method given in \cite{Fo99} for the analogous setting
in the case of the complex Ginibre ensemble, although (\ref{B.6}) is a crucial ingredient made possible by a
recent finding in the literature.

 By definition
 \begin{equation*}
\hat{P}_{N,B}(k) = \left \langle \prod_{l=1}^N {\rm e}^{{\rm i} k b (| \zeta_j|)} \right \rangle,
 \end{equation*}
where the average is with respect to \eqref{Gibbs}. Define
\begin{equation}\label{B.5}
u_l := \int_0^\infty r^{4l+3} {\rm e}^{-2N Q(r)} \, {\rm d}r, \qquad
u_l(b) := \int_0^\infty r^{4l+3} {\rm e}^{-2N Q(r)} {\rm e}^{{\rm i} k b(r)} \, {\rm d}r.
 \end{equation}
According to the theory of radially symmetric skew orthogonal polynomials in the complex plane~\cite[Corollary 3.3]{akemann2021skew}, we have
\begin{equation}\label{B.6}
\hat{P}_{N,B}(k) = \prod_{l=0}^{N-1} u_l(b)/ u_l.
 \end{equation}

 To the integrals in (\ref{B.4}) we now apply Laplace's method of asymptotic analysis.
 For large~$n$,~$L$,~$l$
 such that when divided by $N$ a non-zero limiting value results, this method begins by writing in each integrand
 \begin{equation*}
 r^{4l+3} {\rm e}^{-2N Q(r)} = {\rm e}^{- 2 (n + L + 1) \log (1 + r^2) + (4(L+l) + 3) \log r} = : {\rm e}^{f(r)}
 \end{equation*}
 then expands the integrand about the value of $r$ which maximises the exponent, $r_l$ say.
 An elementary calculation shows that to leading order
 \begin{equation*}
 r_l = \sqrt{\frac{L + l}{n - l}}, \qquad f''(r_l) = - \frac{8 (n - l)^2}{(n+L)}.
 \end{equation*}
 Expanding the integrands in (\ref{B.5}) about this point to second order in the exponent
 shows
 \begin{equation*}
 \hat{P}_{N,B}(k) \sim \prod_{l=1}^N {\rm e}^{{\rm i} k b(r_l)} {\rm e}^{- k^2 (b'(r_l))^2/(2 | f''(r_l)|) } \sim {\rm e}^{{\rm i} k \tilde{\mu}_{N, B}} {\rm e}^{-k^2 \tilde{\sigma}_B^2/2},
 \end{equation*}
 where $ \tilde{\mu}_{N, B}$ is given by the right-hand side of (\ref{B.2}) and
 \begin{equation*}
\tilde{\sigma}_B^2= \frac{1}{4} \int_{r_1}^{r_2} r (b'(r))^2 \,{\rm d}r.
 \end{equation*}
 The right-hand side of this latter expression is equivalent to (\ref{B.4}).
 \end{proof}

 The large $N$ functional form (\ref{B.3}) for the characteristic function of $B$ implies that
 the centred linear statistic $B - \tilde{\mu}_{N,B}$ is a mean zero Gaussian with variance
 given by (\ref{B.4}). The structure of the latter is familiar from the study of the fluctuations
 associated with a linear statistic for the complex Ginibre ensemble; see the recent review \cite[Section 3.5]{Fo22c}.

\subsection*{Acknowledgements} The authors thank Markus Ebke for the help with numerical simulations.
Sung-Soo Byun was supported by Samsung Science and Technology Foundation (SSTF-BA1401-51), by a KIAS Individual Grant (SP083201) via the Center for Mathematical Challenges at Korea Institute for Advanced Study, by the National Research Foundation of Korea (NRF-2019R1A5A1028324), and by the POSCO TJ Park Foundation (POSCO Science Fellowship).
Funding support to Peter Forrester for this research was through the Australian Research Council Discovery Project grant DP210102887.

\pdfbookmark[1]{References}{ref}
\LastPageEnding

\end{document}